\documentclass[preprint,12pt]{elsarticle}
\usepackage{amssymb}
\usepackage{amsmath}
\usepackage{graphicx}
\usepackage{subcaption}
\usepackage{xcolor}
\journal{CMAME}
\usepackage{enumitem} 
\usepackage{hyperref} 
\usepackage{lineno}

\usepackage{scalerel}
\DeclareMathOperator*{\A}{\scalerel*{A}{\sum}}

\usepackage[left=1in,right=1in,]{geometry}

\graphicspath{ {./Figures/} }

\begin{document}

\begin{frontmatter}

\title{An implicit octree-based adaptive Material Point Method}

\author[dur]{Robert E. Bird}
\author[dur]{William M. Coombs\corref{cor1}}\ead{w.m.coombs@durham.ac.uk}
\author[dur]{Charles E. Augarde}
\author[tor,dur]{Giuliano Pretti}
\author[arp,dur]{Ted J. O'Hare}

\cortext[cor1]{Corresponding author}

\address[dur]{Department of Engineering, Durham University, South Road, Durham, DH1 3LE, UK.}
\address[tor]{Department of Civil \& Mineral Engineering, University of Toronto, Toronto, Ontario, Canada.}
\address[arp]{Arup, Central Square, Forth Street, Newcastle upon Tyne, NE1 3PL, UK.}

\begin{abstract}
The Material Point Method provides an effective approach for modelling the large deformations that often arise from contact interactions between rigid structures and surrounding continua. However, solving these problems requires accurate representation of the continuum-structure interface, which necessitates high resolution background mesh and material point discretisations. This requirement, combined with evolving continuum-structure interfaces and the fact that most Material Point Method implementations are dependent on structured meshes, can result in large numerical systems and long run times especially when modelling problems in three-dimensions. Motivated by this issue, this paper provides the first octree-based implicit Material Point Method for efficient solution of large deformation continuum-structure interaction problems. The octree background mesh provides a natural way to automatically adapt both the computational mesh and the material point discretisation based on the position of the interaction between the structure and continuum. The new approach is demonstrated on a number of large deformation benchmark and continuum-structure interaction problems, where up to a $5.5$-times speed up and a consequent $21$-times CO$_2$ saving is achieved when running on a HPC compared to results obtained using a conforming mesh.
\end{abstract}

\begin{highlights}
\item The first octree-based implicit Generalised Interpolation Material Point Method
\item Discrete GIMP basis formulation enables a single non-conforming octree mesh
\item Background mesh and material points automatically adapt based on contact interface
\item Small-cut instability mitigated by extending ghost stabilisation to octree meshes
\item Up to $5.5$-times speed up and $21$-times CO$_2$ saving versus conventional meshes
\end{highlights}


\end{frontmatter}



\section{Introduction}\label{sec:intro}
The Material Point Method (MPM) \cite{sulsky_particle_1994} is a well-established continuum numerical framework for the analysis of problems involving large deformations, for which mesh-based methods are limited by mesh distortion \cite{de_vaucorbeil_chapter_2020,augarde_numerical_2021}. MPM has been successfully applied to a wide range of continuum-rigid structure  interaction problems, for example in metal forming \cite{wikeckowski2004material}, snow-structure interaction \cite{stomakhin2013material}, surgery simulations \cite{ou2025cressim} and soft robots \cite{cochevelou2023differentiable}. In particular, MPM is well suited to modelling geomechanics problems involving granular media where significant deformation and the creation and evolution of free surfaces occur, for example avalanche modelling \cite{gaume2018dynamic} and landslides \cite{andersen2010modelling}. The method has also been successfully applied to geotechnical soil-structure interaction problems, such as numerical Cone Penetration Tests (CPTs) \cite{CECCATO2016440,yost2023addressing,martinelli2022relating}, cable burial ploughs \cite{bird_dynamic_2025}, and drag embedment anchors \cite{birdanchors2026}. In contrast to landslide-type simulations which occur over a short time period, many engineering soil-structure interaction applications occur over much longer time periods, making implicit time integration schemes, where larger time steps can be taken, particularly attractive.

However, the MPM suffers from two well-known issues: the cell-crossing instability and the small-cut problem. The cell-crossing instability arises when material points move between adjacent elements across which the background basis gradients are discontinuous. Detailed discussions and solutions can be found in \cite{steffen2008analysis}, where B-spline basis functions are used, and in \cite{bardenhagen_generalized_2004}, where a particle domain is convolved with the basis functions background grid to form modified basis functions, the resulting method being known as the Generalised Interpolation Material Point Method (GIMPM). A key advantage of the GIMPM is that it can be applied directly to existing Cartesian finite element meshes. The GIMPM was subsequently extended to large deformation implicit analysis by Charlton \emph{et al.} \cite{charlton_implicit_2018}. However the GIMP formulation exacerbates the small-cut problem. The small-cut problem is an ill-conditioning issue, also encountered in unfitted, immersed boundary or fictitious domain methods \cite{coombs2026immersed}, which has been resolved in the MPM through both ghost-penalty \cite{burman2010ghost,coombs_ghost_2023} and aggregation, and extension, techniques \cite{COOMBS2025118012,yamaguchi2021extended}. 

Although implicit MPM implementations are well suited to solid–structure interaction problems, their advantage is offset by the poor scaling of the sparse linear solve with the number of degrees of freedom. This is particularly true for large-deformation problems with non-linear material behaviour, where the system is potentially non-symmetric. The linear solve can be performed either directly, for example by LU factorisation using libraries such as UMFPACK \cite{Davis_2004a}, or iteratively with a Krylov solver such as GMRES \cite{saad1986gmres}. Direct solvers scale poorly for three-dimensional problems: the LU factorisation cost scales at best $\mathcal{O}(N^2)$ \cite{davis2006direct}, whereas Krylov iterative solvers scale at best $\mathcal{O}(kN)$, where $k$ is the number of iterations and $N$ the number of unknowns in the system \cite{saad1986gmres}. Krylov solvers can be accelerated with preconditioners, which reduce the number of iterations required, however this saving must be balanced against the cost of constructing and applying the preconditioner itself. For example a diagonal Jacobi preconditioner is relatively cheap, $\mathcal{O}(N)$, but might be less effective than a ILU type preconditioner, which is both more expensive to compute and apply \cite{saad2003iterative}. Multigrid solvers are alternative choice and are staple in most solver libraries \cite{petsc-efficient,yang2002boomeramg,MueLu} and typically have a scaling of $\mathcal{O}(N)$ \cite{hackbusch2013multi} and can also be used as preconditioners. 

\textit{Octree} meshing strategies are widely used in numerical methods for engineering applications, particularly for adaptive mesh refinement to achieve high solution accuracy and for efficient large-scale parallel scalable computations for finite element problems, see for example  \cite{kolev2021efficient,REINARZ2020107251,burstedde2011p4est}. The benefit of octree meshing is that it enables computational effort to be focused spatially. Elements are refined in regions of the domain that require higher fidelity, whilst coarse elements are retained in less critical areas, giving an analysis that is both efficient and accurate. Octrees have a long history as a hierarchical spatial data structure with their conception presented in the 1980s by Jackins \& Tanimoto \cite{jackins1980oct} and Meagher \cite{meagher1982geometric}. 
Several established octree-based implementations also exist that can serve as a framework for applying the finite element method. Notable examples include deal.II \cite{bangerth2007deal} and DMForest \cite{isaac2015recursive}, the PETSc \cite{petsc-efficient} interface to p4est \cite{burstedde2011p4est}. The octree implementation detailed in this paper is inspired by the work of \cite{tu2005scalable}, which was designed for terascale computation. However, no  claim is made here regarding the scalability of the present implementation.

Adaptive background meshes have previously been investigated within the MPM. In the computer graphics community, Gao \emph{et al.}~\cite{Gao} presented an adaptive GIMP formulation on an octree grid, in which a hierarchy of conforming meshes are used to resolve the hanging nodes. This was initially proposed for an explicit small strain formulation by Tan and Nairn~\cite{tan2002hierarchical}. Although Gao \emph{et al.}~\cite{Gao} identified an implicit treatment as future work, such a formulation has not been reported. In light of recent advances in the stability of implicit GIMPMs, a possibility is that an implicit adaptive GIMPM was, at the time, too unstable to be practical. Separately Zhang \emph{et al.}~\cite{ZHANG2021104097} and Sun \emph{et al.}~\cite{sun2020local} developed adaptive 2D MPMs using B-splines. Ma \emph{et al.}~\cite{ma2006structured} introduced a structured grid-refinement scheme for a 2D GIMPM and Lian \emph{et al.}~\cite{lian2012adaptive} presented a 3D adaptive MPM framework, both with explicit time integration. 

This paper details the first octree-based implicit GIMPM, developed for the efficient solution of soil–structure interaction problems and, more broadly, for large-deformation MPM simulations. A robust implementation requires the formulation of several new numerical components, including a numerical evaluation of the GIMP basis, ghost stabilisation of the GIMPM with hanging nodes, and a systematic mesh-adaptivity strategy. The adaptivity strategy is inspired by goal-oriented error estimation \cite{becker1996feed}, in particular its application to lifting bodies \cite{giani2012anisotropic}, and controls the mesh around the soil–structure interface to achieve a prescribed accuracy of the loads on the structure. The method is validated against analytical, experimental and numerical results for soil–structure interaction problems. The next section outlines the underlying MPM, with the new octree-based GIMPM developed in the sections that follow.

\section{Material Point Method}

The three-dimensional implicit dynamic continuum-structure material point approach adopted in this paper follows that of Bird \emph{et al.} \cite{bird_dynamic_2025,birdanchors2026}, which is based on the quasi-static AMPLE (A Material Point Learning Environment) open source code \cite{coombs2020aample}, with origins in the work of Charlton \emph{et al.} \cite{charlton_implicit_2018} and Coombs \emph{et al.} \cite{coombs2020on}. The following subsections outline key aspects of the approach, with references for more information on each part of the implementation. 

\subsection{Continuum governing equations}
The deformable material is represented by a set of material points that carry information about the properties and state of the material, such as constitutive parameters, mass, deformation history, etc. The governing equations controlling the deformation of the deformable material are assembled and solved on a finite element-like background mesh, which is reset at the end of each time step. However, before the background mesh is reset, the positions and deformation history of the material points are updated, with this updated state being used as the starting point for the next time step. This paper is focused on single phase solid stress analysis where the updated Lagrangian discrete weak form for the deformable body can be expressed as
\begin{equation} \label{eqn:MPMweak} 
   \A_{ \forall p} \Bigl([\nabla_x S_{vp}]^{T}\{\sigma_p\} V_p -[S_{vp}]^{T}\{b\} V_p +[S_{vp}]^{T}\{\dot{v}_p\} m_p \Bigr) - \A_{p\in P_c}\left(\{\delta f_{N,vp} ^{\varphi}\} + \{\delta f_{T,vp} ^{\varphi}\}\right) = \{0\},
\end{equation}
where $\{\cdot\}$ and $[\cdot]$ denote vector and matrix quantities, respectively. The subscript $p$ denotes quantities associated with a material point whereas $v$ denotes quantities associated with the vertices, or nodes, of the background grid. $V_p$ and $m_p$ are the material point volume and mass, $\{\dot{v}_p\}$ is the acceleration of a material point, $\{\sigma_p\}$ is the Cauchy stress at a material point and $\{b\}$ the body force applied over the volume. $P_c$ defines the set of material points in contact with the rigid body. The normal and tangential contact forces acting on the material points are $\{\delta f_{N,vp} ^{\varphi}\}$ and $\{\delta f_{T,vp} ^{\varphi}\}$. $[S_{vp}]$ are the basis functions that link a material point, $p$, to a vertex of the background mesh, $v$. Generalised interpolation basis functions \cite{bardenhagen_generalized_2004,charlton_implicit_2018} are adopted in this paper due to their ability to reduce cell crossing instability whilst remaining efficient for large three dimensional analysis. The values of these basis functions are obtained by performing a volume normalised integral of the conventional linear finite element basis functions over a three dimensional domain associated with each material point (see Section~\ref{sec: GIMP} for details).   

\subsection{Constitutive equations}

This paper adopts large deformation elasto-plastic constitutive equations based on  Charlton \emph{et al.} \cite{charlton_implicit_2018} and Coombs and Augarde \cite{coombs2020aample}. The behaviour of the deformable material points are assumed to be governed by a Hencky material representation, with a linear relationship between elastic logarithmic strains and Kirchhoff stress, $[\tau]$, where the permissible Kirchhoff stress space is bound by a perfect plasticity yield surface. These stress/strain measures are adopted as they recover the conventional infinitesimal strain format of stress return/update algorithms, meaning that these small strain algorithms can be used in large deformation analysis without modification \cite{SouzaNeto}. Details of the yield equations and material parameters are given for each numerical example in Section~\ref{sec: numerical examples}.

\subsection{Rigid body representation and contact approach}\label{sec: Rigid body representation and contact approach}

The structures interacting with the deformable material are assumed to be rigid and represented by a closed triangular faceted mesh. The interaction between the rigid body and the continuum is enforced over each triangular facet using a node-to-surface frictional penalty approach. Within this framework, the vertices of the generalised interpolation material point domains provide the \textit{nodes} and the triangular facets of the rigid body define the contact surface. The normal, $\epsilon_N$, and tangential, $\epsilon_T$, penalty constants are
\begin{equation}
    \epsilon_N = 50E_pA_p^0\quad\text{and}\quad \epsilon_T = 25E_p A_p^0.
\end{equation}
$E_p$ is the Young's modulus of the material point in contact with the rigid body and $A_p^0=(V_p^0)^{2/3}$ is a characteristic area associated with the material point defined by its reference volume $V_p^0$. The values of $50$ and $25$ have units $1$/length and have been determined from the authors' experience as a reasonable balance between contact enforcement and numerical robustness \cite{bird_dynamic_2025}. Friction between the rigid body and the continuum is represented by Coulomb's friction law which captures the stick-slip behaviour in the tangential direction. As seen in (\ref{eqn:MPMweak}), the normal and tangential contact forces provide the coupling between the rigid body and the deformable body represented by material points. These forces therefore influence any unconstrained degrees of freedom of the rigid body (see Bird \emph{et al.} \cite{birdanchors2026} for details).  

\subsection{Coupled material point-rigid body dynamic solution}

The governing equations controlling the deformable material represented by material points, the rigid bodies and their coupling are discretised in time using Newmark's method with the parameters $\gamma=1$ and $\beta=\frac{1}{2}$ \cite{birdanchors2026}. An implicit Newton-Raphson (NR) method combined with a line search is used to solve each time step. Here we employ a simple line search  approach that chooses the lower residual corresponding within the interval $\in[0.1,1]\cdot \Delta t$, where $\Delta t$ is the time step size. Dirichlet (material point displacement and rigid body displacement/rotation) boundary conditions are imposed directly on the background mesh and rigid body using a \textit{strong} approach where the degrees of freedom are constrained out of the linear system and the required value imposed directly. Inhomogeneous Neumann boundary conditions are not required for the problems solved in this paper.  

\section{Octree background mesh}

Octree meshing concentrates computational effort spatially, refining elements only where higher fidelity is required and retaining coarse elements elsewhere. This section describes the background octree mesh and its basis functions, and their extension, through discrete integration, to form the GIMP basis functions linking vertices of the background mesh to each material point.

\subsection{Basis functions for trilinear elements in an octree mesh}
\begin{figure}[h!]
    \centering
    \includegraphics[width=0.7\linewidth]{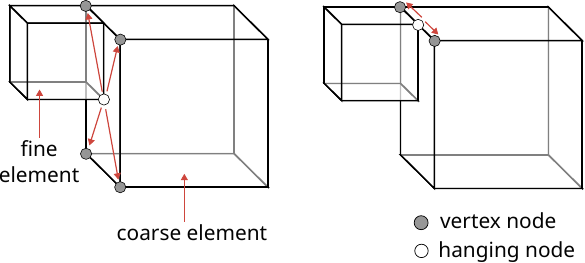}
   \caption{Hanging nodes on a face (left) and on an edge (right). The red arrows link the hanging node to their support nodes.}
    \label{fig:hanging node example}
\end{figure}

The problem domain is the closed cube $\Omega = [0,L]^3 \subset\mathbb{R}^3$, where $L$ is the domain length in each direction, with boundary $\partial\Omega$. The mesh $\mathcal{T}$ is Cartesian of ordered set of hexahedral elements $K\in\mathcal{T}$, with boundaries $\partial K$, such that $\Omega = \bigcup_{K\in\mathcal{T}} K$. To formulate the GIMP basis function it is required that the function space on the background grid is $C^0$ continuous and is formed of elements with a polynomial order of $1$
\begin{equation}\label{equ: polynomial space}
\mathcal{V}(\mathcal{T}) = \left\{ \mathbf{u} \in [C^0({\Omega})]^3 : \mathbf{u}|_K \in [\mathcal{P}_1(K)]^3  \;\; \forall K \in \mathcal{T} \right\}.                               
\end{equation}
where $\mathcal{P}_1(\bullet)$ is the space of degree 1 polynomials on the element $K$. This is not necessarily universal for the MPM, the background grid can be higher order \cite{steffen2008analysis} or a discontinuous space \cite{RENAUD201880}. The background grid elements used here are Lagrange elements, where the degrees of freedom are nodal values at fixed spatial positions and the basis functions satisfy the partition of unity. $C^0$ continuity across element interfaces is satisfied when adjacent elements share nodes, this ensures that data from the nodes is interpolated from adjacent elements to the same values on the shared face. An octree mesh introduces \textit{hanging} nodes, nodes that lie on an face, or edge, of a coarser neighbouring element but do not coincide with its vertices, breaking the continuity of the approximation across these interfaces and creating a discontinuity and therefore no longer satisfying Equation~\eqref{equ: polynomial space}. Methods to resolve this are well established in the literature and have been well addressed by many authors with constraint equations \cite{solin2010adaptive,bangerth2007deal}, the use of B-splines also applied to MPM \cite{sun2020local,ZHANG2021104097}, or Nitsche's method in discontinuous Galerkin methods \cite{arnold2002unified}. Here we use the constraint equation approach as it can be naturally included in the GIMP formulation directly.  First, it is necessary to define when a node is considered hanging. Let $\mathcal{N}(K)$ be the set of vertex nodes of element $K$, and $\mathcal{N}(\mathcal{T})$ the set of all nodes in the mesh. A node $n$, with position $\mathbf{x}_n$, is hanging if its position is on the boundary of an element but it is not a vertex node of that element,
\begin{equation}\label{equ: hanging node equation}
     \mathbf{x}_{n} \in \partial K \quad \text{and} \quad n \notin \mathcal{N}(K).
\end{equation}
The set of hanging nodes on the element $K$ is therefore defined
\begin{equation}
        \mathcal{N}_H(K) = \{n\in \mathcal{N}(\mathcal{T}): \mathbf{x}_{n} \in \partial K ~~ \text{and} ~~ n \notin \mathcal{N}(K)\}
\end{equation}
with the set of all (unique) hanging nodes and non-hanging nodes in the mesh
\begin{equation}
    \mathcal{N}_H(\mathcal{T}) = \bigcup_{K\in\mathcal{T}} \mathcal{N}_H(K)\quad\text{and}\quad  \mathcal{N}_V(\mathcal{T}) = \mathcal{N}(\mathcal{T})\setminus \mathcal{N}_H(\mathcal{T}).
\end{equation}
A node $n$ can be a hanging node on a element face or edge, see Figure \ref{fig:hanging node example} as an example.
The coarse element on which the node $n$ is hanging is defined $K_n$, for a face hanging node the mapping between hanging node $n$ and $K_n$ is unique, however this is not the case for an edge hanging node where it is possible that there are three elements. To make the mapping from the node $n$ to $K_n$ unique the following definition is required
\begin{equation}
    K_n = \max_{K\in\mathcal{T}}\{K:\ n\in \mathcal{N}_H(K)\}.
\end{equation}
For an 8-noded hexahedral element, the approximation of $\mathbf{u}$ at $\mathbf{x} \in K$ is
\begin{equation}
    \mathbf{u}(\mathbf{x}) = \sum_{n\in\mathcal{N}(K)} \phi_{nK}(\mathbf{x}) \, \mathbf{u}_n,
\end{equation}
where $\phi_{nK}$ are the trilinear basis functions and $\mathbf{u}_n$ the nodal values for the element $K$. If a node $n$ on the fine element $K$ is identified as hanging, $n\in\mathcal{N}_H(\mathcal{T})$, the nodal value is not independent but is constrained to match the approximation of the coarse neighbouring element $K_n$,
\begin{equation}\label{equ: constraint equation}
    \mathbf{u}_{n} = \sum_{i \in \mathcal{N}(K_n)} \phi_{iK_n}(\mathbf{x}_{n}) \, \mathbf{u}_{i},
\end{equation}
where $\phi_{iK_n}$ are the basis functions of $K_n$ evaluated at the hanging node position $\mathbf{x}_{n}$, and $\mathbf{u}_i$ are its nodal values, see Figure \ref{fig:hanging node example}. Equation~\eqref{equ: constraint equation} is general and for trilinear elements can be reduced to consider only the basis functions which are non-zero at $\mathbf{x}_{n}$,
\begin{equation}\label{equ: constraint equation reduced}
    \mathbf{u}_{n} = \sum_{i \in \mathcal{S}(n)} \phi_{iK_n}(\mathbf{x}_{n}) \, \mathbf{u}_{i},
\end{equation}
where
\begin{equation}
\mathcal{S}(n) = \{ i \in \mathcal{N}(K_n) : \phi_i(\mathbf{x}_{n}) \neq 0 \}
\end{equation}
returns the vertex nodes which support the hanging node and
\begin{equation}
\mathcal{Q}(i,K) = \{n\in\mathcal{N}_H(\mathcal{T}): i\in\mathcal{S}(n)\}\cap\mathcal{N}(K)
\end{equation}
is the set of hanging nodes which $i$ supports which are also vertices to $K$. The set of all independent nodes in an element is the union of the nodes which are only vertices, $\mathcal{N}_V(K) =\mathcal{N}(K) \setminus \mathcal{N}_H(\mathcal{T})$ and the nodes which support the hanging nodes,
\begin{equation}
\mathcal{I}(K) =\mathcal{N}_V(K) \cup\left(\bigcup_{n \in \mathcal{N}(K)\setminus \mathcal{N}_V(K)} \mathcal{S}(n) \right)
\end{equation}
see Figure \ref{fig:hanging node example} for examples of the vertex support nodes for edge and face hanging nodes. The modified basis for the active node $i\in\mathcal{I}(K)$ is
\begin{equation}\label{equation: hanging node basis}
    \tilde{\phi}_{iK}(\mathbf{x}) = \phi_{iK}(\mathbf{x}) + \sum_{\substack{n\,\in\,\mathcal{Q}(i,K)}}
    \phi_{nK}(\mathbf{x})\, \phi_{iK_n}(\mathbf{x}_n)
\end{equation}
This allows the solution $\mathbf{u}$ in $K$ to be simply written as
\begin{equation}\label{equation: hanging node solution}
    \mathbf{u}(\mathbf{x})
    = \sum_{i \in \mathcal{I}(K)}
    \tilde{\phi}_{iK}(\mathbf{x}) \mathbf{u}_i.
\end{equation}

\subsection{Generalised interpolation material point basis functions}\label{sec: GIMP}
In the standard MPM the material point basis functions are the background grid basis. However for a $C^0$ basis, when a material point crosses the boundary of adjacent elements the derivative of the basis for the shared nodes on the interface will jump in their value, leading to the so-called cell-crossing instability. A review of the different methods to resolve this is presented in the introduction. Here the GIMPM \cite{bardenhagen_generalized_2004,charlton_implicit_2018} is used, where the basis functions are constructed by convolution of a particle domain with a $C^0$ background grid basis to form a new basis which is $C^1$.

The general calculation of a GIMP basis function with a unity characteristic function over the domain associated with the material point is
\begin{equation}\label{equ: GIMP basis}
S_{vp} = \frac{1}{V_p} \sum_{K\in\mathcal{K}_p}\int_{\Omega_p \cap \Omega_K} \tilde{\phi}_{vK}(\mathbf{x}) \, d\mathbf{x}\quad\text{where}\quad  \Omega_K(v) = \bigcup \left\{ \Omega_K : v \in \mathcal{I}(K) \right\},
\end{equation}
where
\[
\mathcal{K}_p = \{ K \in \mathcal{T} : \Omega_K \cap \Omega_p \neq \emptyset \},
\]
and
$V_p$ is the volume of the domain associated with the GIMP $p$, $v$ is the global node number, $\Omega_p$ is the domain of $p$, $\Omega_K$ is the domain of the element $K$, $\Omega_K(v)$ is the set of element domains which contain the node $v$ and $\tilde{\phi}_{vK}$ is the basis function for the node $v$ in element $K$. 

\begin{figure}[ht!]
\centering
\begin{subfigure}[b]{0.4\textwidth}
      \centering
          \raisebox{1.5cm}{\includegraphics[width=0.97\textwidth]{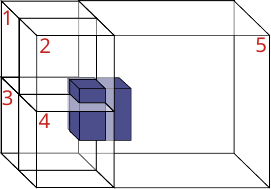}}
          \caption{}
          \label{fig:gimp_part1}
          \end{subfigure}
      \hfill
      \begin{subfigure}[b]{0.58\textwidth}
      \centering
          \includegraphics[width=\textwidth]{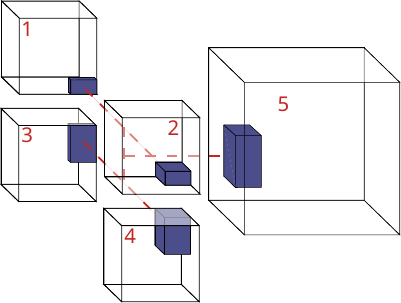}
          \caption{}
          \label{fig:gimp_part2}
      \end{subfigure}
      \caption{GIMP basis functions: (a) is shows the intersection of a GIMP domain (blue) with the background mesh, (b) is an exploded view of the GIMP intersections with each of the five elements it intersects.}
      \label{fig:gimp_integration}
\end{figure}  

For conforming meshes the GIMP basis is normally presented as a three-part piecewise function of a 1D GIMP moving over element intervals. A tensor product along the Cartesian directions is then used to create the trilinear basis in 3D, see for example Charlton \emph{et al.} \cite{charlton_implicit_2018}. For the non-conforming mesh of the octree, building a piecewise function for all possible variations of GIMP intersections with the background mesh, and the subsequent GIMP basis function, would be tedious and prone to error in terms of presentation, interpretation and implementation. Instead, since the background grid basis, $\tilde{\phi}_{vK}$, is linear inside the domain of the element $K$, a single Gauss point quadrature rule can be used to compute the $S_{vp}$. Where the numerical integration of Equation~\eqref{equ: GIMP basis} is
\begin{equation}\label{equ: GIMP numerical}
S_{vp} = \frac{1}{V_p} \sum_{K\in(\Omega_p \cap \Omega_K(v))} \tilde{\phi}_{vK}(\mathbf{x}_{Kp})V_{Kp}
\end{equation}
where $\mathbf{x}_{Kp}$ is the centre of the intersection $\Omega_p \cap \Omega_K(v)$, see Figure \ref{fig:gimp_part2}, and $V_{Kp}$ is the volume, $|\Omega_p \cap \Omega_K(v)|$.

\begin{figure}[ht!]
\centering
\begin{subfigure}[b]{0.55\textwidth}
      \centering
          \raisebox{0.45cm}{\includegraphics[width=\textwidth]{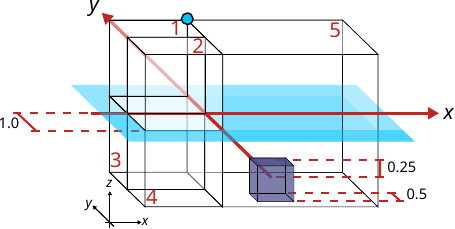}}
          \caption{}
          \label{fig:gimp_example_mesh}
          \end{subfigure}
      \hfill
      \begin{subfigure}[b]{0.44\textwidth}
      \centering
          \includegraphics[width=\textwidth]{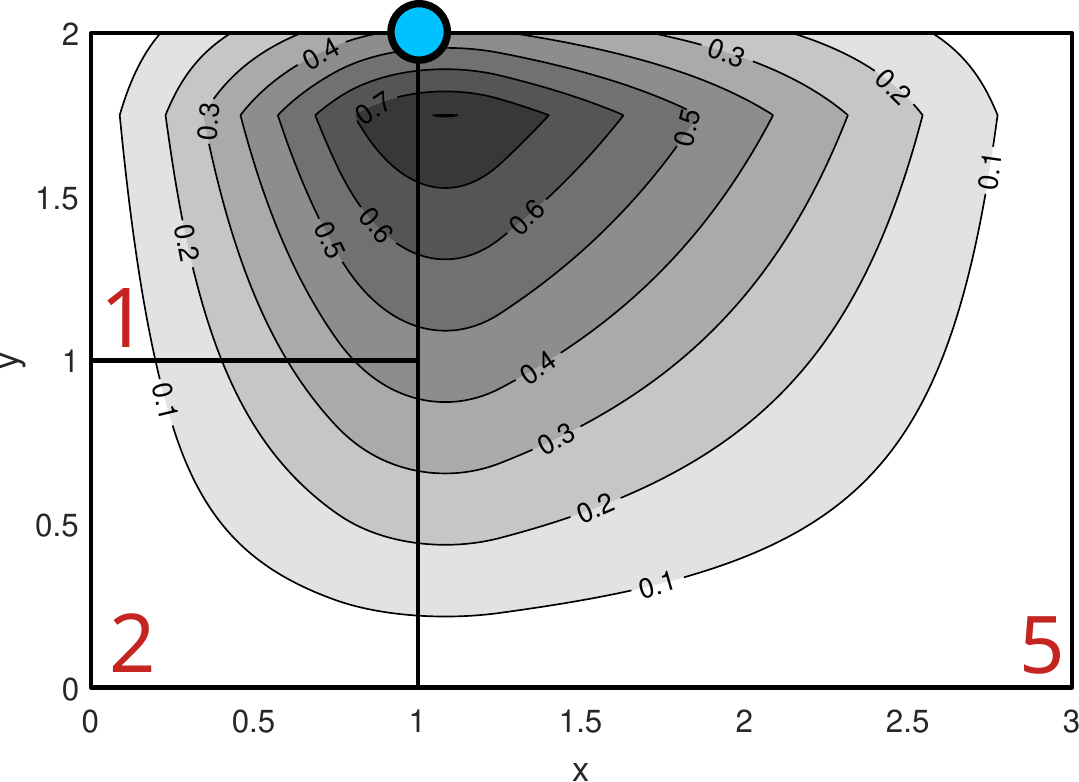}
          \caption{}
          \label{fig:gimp_example_part_results}
      \end{subfigure}
      \caption{GIMP: (a) is a schematic of five-element mesh and a GIMP moving over the plane marked in blue, (b) is the value of $S_{vp}$ over the blue domain and vertex node of (a).}
      \label{fig:gimp_example_mesh and gimp_example_part_results}
\end{figure}

To demonstrate how the GIMP basis function $S_{vp}$ varies for a vertex node $v$ which supports a hanging node, a GIMP, of side length $0.5$, is moved over the $x$-$y$ plane which intersects the vertex node, marked in blue, and the hanging node, marked in grey, Figure \ref{fig:gimp_example_mesh}. The GIMP is not varied in $z$ so that its centre point always coincides with the plane. Figure \ref{fig:gimp_example_mesh and gimp_example_part_results}(b) is a plot of the blue $x-y$ plane, it shows how the value of $S_{vp}$ varies for the vertex node, the blue node in Figure \ref{fig:gimp_example_mesh}, as the GIMP is moved on the plane. Figure \ref{fig:gimp_example_mesh and gimp_example_part_results}(b) plots $S_{vp}$ on the blue $x$–$y$ plane, showing how its value at the vertex node (the blue node in Figure \ref{fig:gimp_example_mesh}) varies as the GIMP is moved across the plane. The plot demonstrates how constraining hanging nodes causes the vertex node to influence elements with which it does not coincide, and also how the $C^1$ continuity is also maintained.

\section{Ghost stabilisation}\label{sec: ghost}
As introduced in Section~\ref{sec:intro}, the small cut instability is a problem in unfitted or fictitious domain finite element methods and occurs where a domain boundary cuts through an element, leaving a small portion inside the problem domain, resulting in very small contributions to the mass and or stiffness matrix and the loss of coercivity of the linear system of equations. Coombs \cite{COOMBS2025118012} extended the so-called ghost stabilisation solution in Burman \cite{burman2010ghost} to MPM and GIMP for implicit and dynamic-explicit solution methods. Other stabilisation methods exist for MPMs, for example aggregation \cite{COOMBS2025118012}, which works by constraining degrees of freedom on the problem boundary. However, the approach used for considering hanging nodes already constrains the solution at hanging nodes to the their neighbouring vertex nodes; applying a second constraint would make the framework cumbersome and hence a more natural method to use here is ghost stabilisation.
\begin{figure}[ht!]
\centering
\begin{subfigure}[b]{0.3\textwidth}
      \centering
          \includegraphics[width=\textwidth]{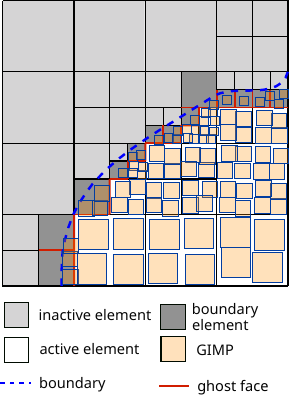}
           \caption{}
          \label{fig:2D hanging node mesh}
          \end{subfigure}
      \hfill
      \begin{subfigure}[b]{0.68\textwidth}
      \centering
          \raisebox{0.5cm}{\includegraphics[width=\textwidth]{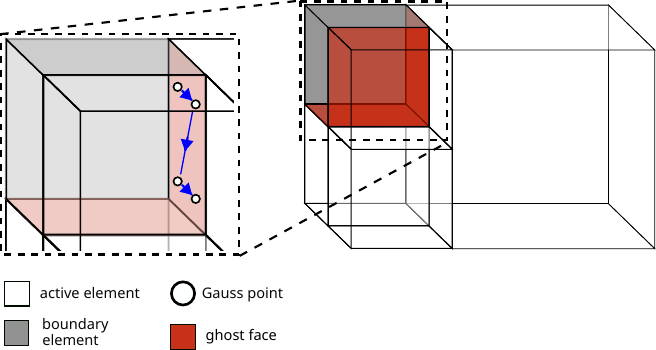}}
          \caption{}
          \label{fig:3D hanging node mesh}
      \end{subfigure}
      \caption{Ghost stabilisation: (a) is a 2D mesh with hanging nodes showcasing where the ghost stabilisation is applied with (b) extending this example to 3D but for a simpler mesh with zoomed view of a face being integrated with Gauss quadrature.}
      \label{fig:2D hanging node mesh and 3D hanging node mesh}
\end{figure} 
Consider a mesh $\mathcal{T}$ with elements $K \in \mathcal{T}$. Let $\mathcal{F}_I$ denote the set of interior faces of the mesh, where each face $F \in \mathcal{F}_I$ is defined as $F = K^+ \cap K^-$ for two neighbouring elements $K^+, K^- \in \mathcal{T}$, values associated with the elements are respectively denoted with the superscripts $(\bullet)^+$ and $(\bullet)^-$. For a given element $K$, we denote by $\mathcal{F}_K \subset \mathcal{F}_I$ the subset of interior faces belonging to $K$.  Here the ghost stabilisation is only applied to the stiffness matrix, but it can also be applied to the mass matrix, see Sticko \emph{et al.} \cite{sticko2020high}. The ghost bilinear form for the stiffness matrix stabilisation of element $K$ is,
\begin{equation}\label{eq:ghost}
    j_F(\boldsymbol{u}_h, \boldsymbol{v}_h)=\frac{h_F\beta_F}{\kappa} \left( [\partial^k_n \mathbf{u}],\, [\partial^k_n \mathbf{v}] \right)_F
\end{equation}
where 
\[
 [\partial^k_n (\bullet)] = \left(\frac{\partial(\bullet)^+}{\partial \mathbf{x}} - \frac{\partial(\bullet)^-}{\partial \mathbf{x}}\right)\cdot\mathbf{n}^+,
\]
$\boldsymbol{v}_h$ and $\boldsymbol{u}_h$ are the test and trial functions, $\kappa = 10$ is a constant, and $F$ is element of the set of ghost faces $\mathcal{F}_G\subseteq\mathcal{F}_I$. $h_F = \max(h_{K^+}, h_{K^-})$ is the maximum element size with $h_K$ the maximum side length of element $K$, a similar approach to \cite{massing2014stabilized} where the average element size is used, and $\beta_F = \max(\beta_{K^+}, \beta_{K^-})$ is the maximum Young's modulus where $\beta_K$ the maximum Young's modulus of the material points within $K$.

To obtain the ghost faces for the GIMP with a mesh with hanging nodes the same algorithm as Coombs \cite{coombs_ghost_2023} can be used but applied to here to an octree mesh, see Figure \ref{fig:2D hanging node mesh and 3D hanging node mesh} for a 2D example of mesh with the small-cut problem. Let $\mathcal{T}^{\emptyset} \subset \mathcal{T}$ and $\mathcal{T}^p \subset \mathcal{T}$ denote the sets of unpopulated and populated elements respectively. The set of boundary elements, $\mathcal{T}^B\subset\mathcal{T}^p$ are defined as those sharing a face with an unpopulated element,
\begin{equation}
    \mathcal{T}^B = \{ K \in \mathcal{T}^p : \text{there exists}\, F \in \mathcal{F}_K,\ F = K \cap K',\ K' \in \mathcal{T}^{\emptyset} \},
\end{equation}
and the ghost stabilisation faces are the faces of boundary elements shared with either another boundary element or a populated element,
\begin{equation}
    \mathcal{F}_B = \{ F = K^+ \cap K^- : K^+ \in \mathcal{T}^B,\ K^- \in \mathcal{T}^p \}.
\end{equation}
Care is required when integrating over a ghost face if the elements sharing the face $F$ are of different size, see Figure~\ref{fig:3D hanging node mesh}, marked with the white circles indicating Gauss quadrature. As the integration over the face of the smaller element is over its entire side face, the Gauss quadrature is always defined for the smaller face and then mapped to the subregion of the larger element. This way the Gauss quadrature can be defined for the faces of an element at the beginning of the algorithm, and a mapping procedure calculates the position of the Gauss points on the larger neighbour element.

\section{Adaptive strategy}\label{sec:adaptive stratergies}
\begin{figure}[ht!]
    \centering
    \includegraphics[width=\textwidth]{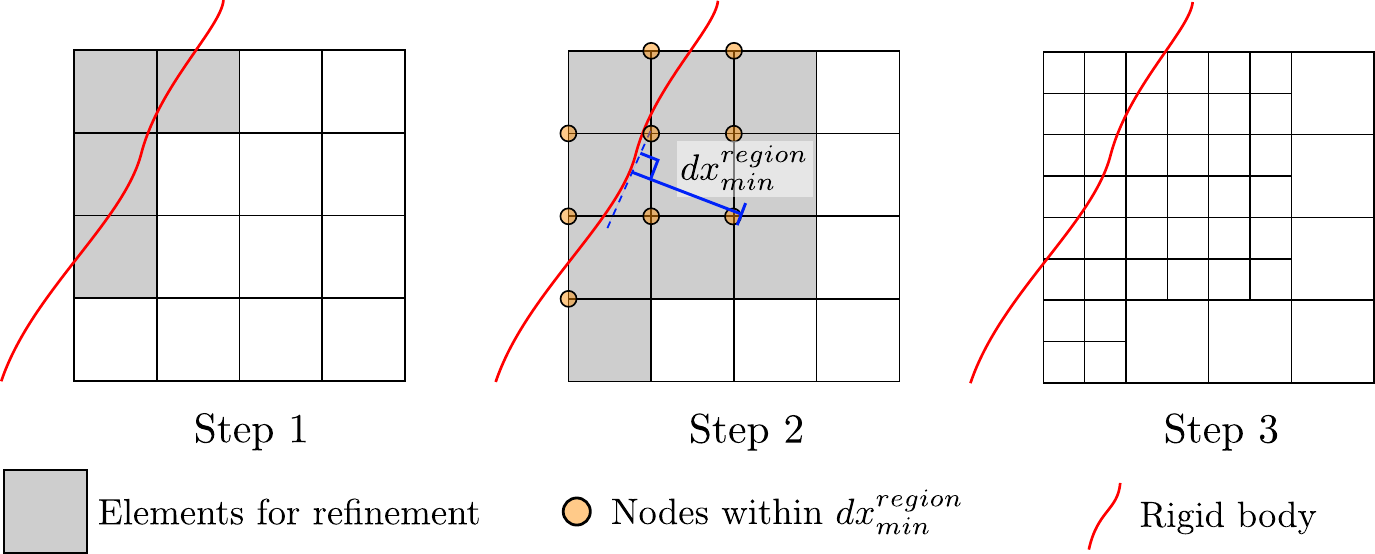}
    \caption{Adaptive strategy: The three steps of element adaptivity in the proximity of a rigid body.}
    \label{fig:ref scheme}
\end{figure}

At the start of each load or time step the octree mesh is rebuilt with an adaptivity process that occurs in two stages: refining the mesh and then the GIMPs. Rebuilding the mesh each time means it can be adaptively refined to match the current state directly, this removes the need for a separate derefinement algorithm to coarsen regions that no longer require fine resolution. It also fits into the framework of the Material Point Method; all the Lagrangian data is stored at Lagrangian points whilst the background mesh is reset at the beginning of each step. Here the adaptivity scheme is governed by the position of where the rigid body surface interacts with the GIMPs, as shown in Figure \ref{fig:ref scheme}. This approach is inspired by goal-oriented error estimators, for example Giani \cite{giani2012anisotropic}, where mesh refinement is localised, as indicated by the error estimator, around the interface of the body to maximise the accuracy of the computed load on a solid object, such as the lift and drag on an aerofoil. This is a deliberate choice motivated by the focus on the accurate evaluation of loads on the rigid bodies from deformable–rigid body interaction. It is noted that this marking scheme could be replaced by a different metric for other analyses. The process for refining elements occurs in the three steps outlined in Figure \ref{fig:ref scheme}. 
\begin{itemize}[leftmargin=2cm]
      \item [Step 1:] Identify the elements that intersect the rigid body and mark these for refinement;
      \item [Step 2:] Find all nodes that are within the distance $dx_{min}^{region}$. Any elements with these nodes are marked for refinement;
      \item [Step 3:] Refine all elements which have been identified for refinement.
\end{itemize}
During Step 2 a choice is made here that no more than one hanging node can exist on an edge or face, to keep the adaptive numerical framework simple, \cite{solin2003higher}, but which is not a numerical restriction \cite{solin2010adaptive}. To obtain a more refined mesh, Steps 1-to-3 are repeated until the desired element size is achieved. 

\begin{figure}[ht!]
    \centering
    \includegraphics[width=0.9\linewidth]{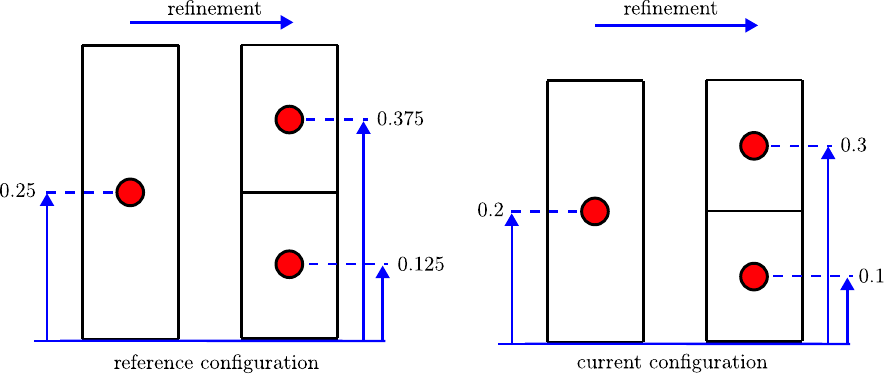}
    \caption{Adaptive strategy: One dimensional refinement example, reference and current configurations.}
    \label{fig:1D ref example}
\end{figure}

\begin{figure}[ht!]
    \centering
    \includegraphics[width=0.6\textwidth]{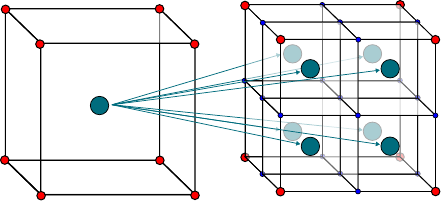}
    \caption{Adaptive strategy: How material information which is constant throughout the GIMP domain is projected from the coarse to refined GIMPs.}
    \label{fig:GIMP_refinement_centre}
\end{figure}

The GIMP refinement is defined by the sizes of the elements it occupies. The set of elements that intersect the GIMP $p$ is defined as
\[
\mathcal{K}_p = \{ K \in \mathcal{T} : \Omega_K \cap \Omega_p \neq \emptyset \},
\]
where the \textit{half} side length of the GIMP is $l_p$ and the half side length of the element $K\in\mathcal{K}_p$ is  $l_K$. A GIMP is refined if
\begin{equation}
    \exists K\in\mathcal{K}_p:l_K<\beta ~l_p
\end{equation}
where $\beta\geq2$ and is user defined (for all problems considered here $\beta=2$). Once a GIMP is refined it is necessary to transfer its data to the new material points, as shown in Figure \ref{fig:GIMP_refinement_centre}. The material and kinematic data associated with a GIMP is constant within its domain therefore, when a larger GIMP is split into 8 smaller GIMPs the data are directly transferred to the eight new GIMPs. However, when the material data depends on the position in the reference configuration, each point generated by refinement must evaluate its own material data from its reference coordinates, so that it is initialised consistently with the points present in the original discretisation. Consider a 1D example shown in Figure \ref{fig:1D ref example} where the Young's modulus, $E$, has a linear variation with the vertical position $E = (1-z)$, where $z$ is the vertical position. 
The vertical position of a GIMP is defined by its centre point, the domain length of the GIMP in the vertical direction is referred to as the GIMP length. In the current configuration the GIMP has vertical position $0.2$ and length $0.4$, in the reference configuration these are $0.25$ and $0.5$. Refinement generates two new material points which, in the current configuration, have a length $0.2$ and are at vertical positions $0.1$ and $0.3$. These still represent a coarse discretisation of $E$ and must be updated to remain consistent with the criterion used to assign $E$, therefore their equivalent vertical positions in reference configuration are used, $0.125$ and $0.375$. This update warrants further investigation. In particular, it is not derived from an energetic framework, which is problematic for plastic materials as it introduces an inconsistency between the updated material properties and the kinematic variables. Nevertheless, the numerical results presented here suggest this inconsistency is small, as good agreement is obtained with previously published results.

Additionally, no transfer of contact data is necessary. This is because the most refined mesh is always associated with the elements in the region around the rigid body and hence the material points that could be in contact are never refined further. Finally, although the mesh is derefined as the rigid body moves through it, the GIMPs themselves are not. Derefining the GIMPs would require approximating the Lagrangian data of several small GIMPs by lumping them back into a single larger GIMP. This was considered unnecessary, since computing the residual and tangent matrix of a material point is cheaper and has the potential for a higher degree of parallel processing, compared with the linear solve when the number of degrees of freedom ($N$) becomes large. For example a direct solve costs scales with $\mathcal{O}(N^2)$ \cite{davis2006direct}, and iterative solvers are at best $\mathcal{O}(kN)$, where $k$ is the number of iterations \cite{saad1986gmres}. Whereas the cost of computing GIMPs is relatively small, Feng \emph{et al.} \cite{feng2026mpm} reported that, although the construction cost is non-negligible, approximately $64\%$ of the computational time in their implementation was spent in the linear solve.

\section{Numerical examples} \label{sec: numerical examples}

In this section, five numerical examples are presented. The first two examples, the column under self weight and the column under constant stress from contact, validate the GIMP shape formulation of Equations \eqref{equation: hanging node basis} and \eqref{equ: GIMP basis}. The meshes will have hanging nodes to demonstrate convergence with mesh and GIMP refinement when compared to the analytical solutions. The third example is an implicit-in-time dynamic problem of a ball rolling down a slope where the position of the sphere at time $t$ is compared to the analytical solution for different friction coefficients, also inspected with explicit-MPM by \cite{su2025pinball}. The analytical solution considers two rigid bodies, however here the slope is made sufficiently stiff to approximate it as a rigid body to validate the hanging node formulation with contact and friction. The last two problems are validations against experimental data: a Cone Penetration Test (CPT) and a drag anchor pull. The numerical results are compared to experimental data whilst the simulation runtimes are compared to the results of equivalent accuracy but run with conforming meshes, data  originally presented in Bird \emph{et al.} \cite{bird_dynamic_2025,birdanchors2026}.

For timing results it is important to note that all simulations were performed on a single compute node of Hamilton8, provided by Hamilton HPC Service of Durham University, equipped with two AMD EPYC 7702 64-core processors (128 cores in total) and 256 GB of RAM. Each simulation was allocated 10 cores and up to 80 GB of RAM. For green energy calculations it was assumed that the full allocation of resources was always used.

\subsection{One dimensional column compression under self weight}
\begin{figure}[ht!]
    \centering
    \includegraphics[width=0.5\linewidth]{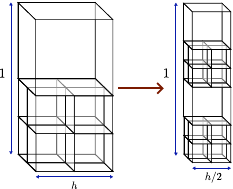}
   \caption{Compression under self weight: An example of the refinement scheme with hanging nodes.}
    \label{fig:example mesh}
\end{figure}

Convergence of the implicit GIMPM (iGIMPM) problem was first demonstrated by Charlton \emph{et al.} \cite{charlton_implicit_2018} on the 1D self-weight column. In that paper it was shown that the error in the stress solution, normalised with respect to the volume, converges with uniform refinement for a conforming mesh. Here the same general problem and GIMP domain update is considered as in Charlton \emph{et al.} \cite{charlton_implicit_2018}. The initial mesh, and a refinement step, are shown in Figure \ref{fig:example mesh}. 
The domain has dimensions $(x,y,z)\in[0,W]\times[0,W]\times[0,0.8]$ m, where $W$ is the width of the domain and is equal to the side length of the largest element in the domain. The first mesh has $W=0.4$ m, and is divided by 2 each refinement step, as in Figure \ref{fig:example mesh}. When reporting on the convergence of the problem the \textit{size} of elements in the mesh is reported as $W$. On all sides of the domain, roller boundary conditions exist except the top where there is a homogeneous Neumann boundary condition, all nodes are at least fixed in $x$ and $y$ to make the problem one-dimensional. The validation uses slightly different material properties to that of Charlton \emph{et al.} \cite{charlton_implicit_2018} so that significant deformation occurs and overlap of GIMPs between elements of different sizes occurs. The material is Hencky elastic, with constant elastic parameters, Young's modulus $E=10^3$ Pa, Poisson's ratio $\nu = 0.0$, and the material density is set to $\rho = 50$ kg/m$^3$. The domain is subject to a body force from gravitational loading $g_i = [0,0,-9.81]$ m/s$^2$, and the problem is solved using a Newton-Raphson scheme, where the load is incrementally increased over 20 load steps. 

To demonstrate convergence, the numerical stress solution is compared to the analytical solution and a stress error is defined which is normalised to the domain volume. The measure is similar to that of \cite{charlton_implicit_2018}, but normalised with respect to the stress at the bottom of the column, i.e.
\begin{equation}
    e_{\sigma} = \frac{1}{\sigma_gV_\Omega}\sum_{p\in P} |\sigma^z(z^0_p)-\sigma_p^z|V_p
\end{equation}
where $\sigma_g = 0.8\rho g$ is the stress at the bottom of the domain, $V_\Omega$ is the domain volume, $\sigma^z(z_p^0) = \rho g(0.8-z_p)$ is the vertical stress solution $z_p^0$ is the vertical position of the point at time $t=0$ and $V_p$ is the GIMP volume. For this problem in the smaller elements the GIMPs are initiated in a $2\times2\times2$ structure, whereas in the larger elements the structure is $4\times4\times4$. This is to ensure no GIMP refinement occurs during the validation, focusing the study on the convergence of the method with an octree background mesh with hanging nodes.

\begin{figure}[ht!]
    \centering
    \includegraphics[width=0.7\textwidth]{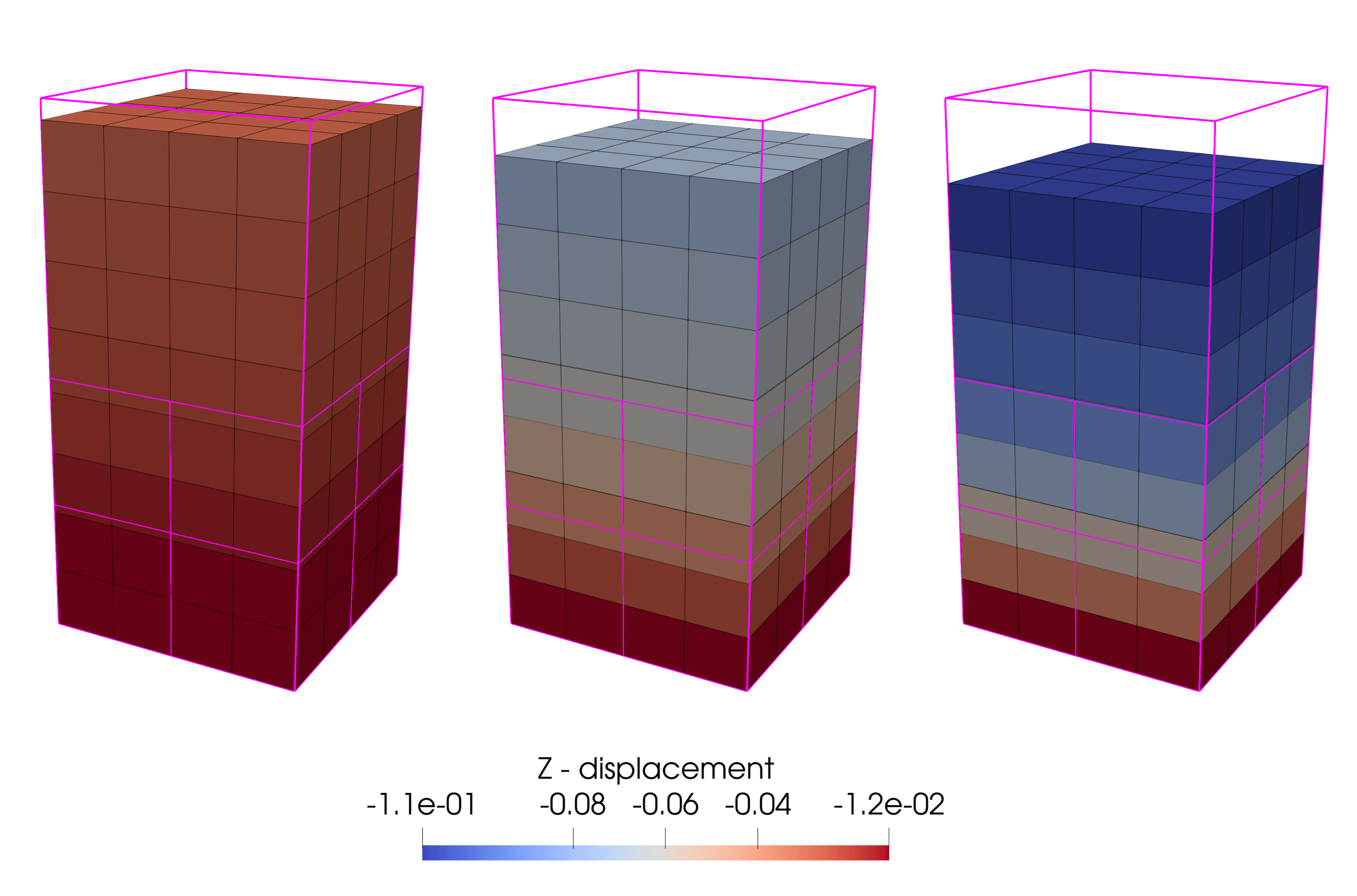}
    \caption{Self weight: Displacement plot of the GIMPs and mesh for steps 1, 10 and 20, of 20.}
    \label{fig:column_stress_solution}
\end{figure}

An example of a column, $W=0.4$ m, undergoing deformation over the 20 load steps is shown in Figure \ref{fig:column_stress_solution}. The mesh is the pink wire frame and the GIMPs are shown by the blue-to-red coloured cuboids, with blue being the most negative vertical displacement. The figures shows that the GIMPs are traversing the element boundaries, particularly at $z=0.4$ m, where hanging nodes exist.

\begin{figure}[ht!]
    \centering
    \includegraphics[width=0.5\textwidth]{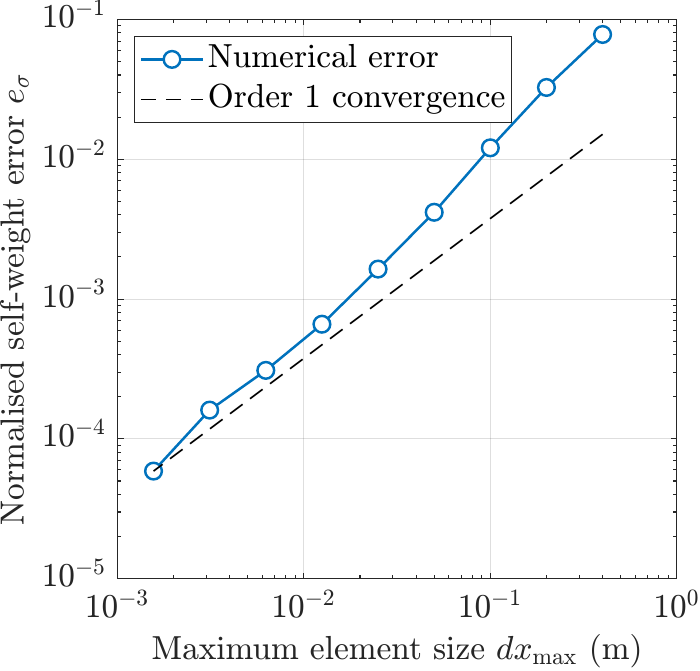}
    \caption{Self weight: Convergence of the error with mesh refinement.}
    \label{fig:column_convergence}
\end{figure}

Convergence of the solution with mesh refinement is shown in Figure~\ref{fig:column_convergence}, plotted against a reference slope of one. As in Charlton \emph{et al.}~\cite{charlton_implicit_2018}, the observed convergence 
rate of the stress solution is slightly greater than unity, whereas linear finite elements should yield exactly first-order convergence. This is caused by the GIMP basis functions having different orders when the domain of a material point is fully within an element (linear) or overlapping multiple elements (quadratic) \cite{coombs2020on}. 

\subsection{Column under constant stress from contact}\label{sec: contact problem}
The purpose of this problem is to validate the GIMPM formulation for normal contact between a rigid and deformable body, originally presented by Bird \emph{et al.} \cite{bird_dynamic_2025}, and apply it here to problems where the GIMP contact is occurring in elements with hanging nodes. The validation is for penalty based contact and will show: 
\begin{enumerate}
    \item as the penalty is increased the overlap between the bodies will converge to zero;
    \item the stress solution in the deformable bodies converges to the analytical solution; and
    \item the vertical stress is the constant throughout the deformable contact plane and the rest of the body.
\end{enumerate}

\begin{figure}[ht!]
\centering
\begin{subfigure}[b]{0.49\textwidth}
      \centering
          \raisebox{0.0cm}{\includegraphics[width=0.7\textwidth]{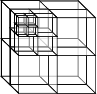}}
           \caption{}
          \label{fig:column_compression_series}
          \end{subfigure}
      \hfill
      \begin{subfigure}[b]{0.49\textwidth}
      \centering
          \includegraphics[width=0.7\textwidth]{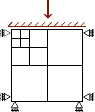}
          \caption{}
          \label{fig:column_compression_series_BCs}
      \end{subfigure}
      \caption{Contact cube: The initial mesh is shown in (a) with a 2D representation of the boundary conditions and application of the rigid body shown in (b).}
      \label{fig:gimp_example_part}
\end{figure}

The material is Hencky elastic with a Young's modulus of $E=10^6$ Pa, and Poisson's ratio of $0.0$. All degrees of freedom are fixed in the $x$ and $y$ directions and roller boundary conditions exist on all faces, except the top face which is homogeneous Neumann. The mesh is kept constant throughout this study and has dimensions $(x,y,z)\in[0,0.8]^3$ m, see Figure \ref{fig:column_compression_series}, and is initiated with $2\times2\times2$ material points within each element. The rigid body contacts the problem from the top of the domain and impinges on it a distance $0.2$ m, marked by the red line in Figure \ref{fig:column_compression_series_BCs}, applied over 20 load steps.

The domain should undergo a uniform compression and therefore the error $e_\sigma$ can be computed by comparing the analytical stress solution to the numerical solution
\begin{equation}\label{eq: L2 stress error}
    e_{\sigma} = \left({ \sum_{\forall p} |\sigma^z - \sigma_{p}^z|^2 V^0_{p} } \right)^{1/2} \quad\text{with}\quad \sigma^z = E \left(\frac{l_0}{l}\right)\log\left({l}/{l_0}\right)
\end{equation}
where $\sigma^z$ is the analytical stress solution in the vertical direction, $\sigma_p^z$ is the numerical vertical stress for the point $p$, $V_p^0$ is the original GIMP volume, $l_0$ is the original length and $l$ is the final length. The contact penalty in the normal direction is defined as
\[
 \epsilon_N = p_fE_pA_p^0
\]
where $E_p$ is the Young's modulus of the GIMP in contact, $A_p^0 = (V_p^0)^{2/3}$ is the contact surface area of the point, and $p_f$ is the penalty factor, the parameter which is varied to study the convergence of the problem. The second error measure used is the position error, which measures the difference in the final minimum $z$-position of the rigid body and compares it with the final position of the top face of the GIMPs in contact, i.e.
\begin{equation}
     e_u = \left(\sum_{p\in P_c} (u_{p} - u_{RB})^2A_p)\right)^{1/2}
\end{equation}
where $P_c$ is the set of material points $p$ in contact, $u_p$ is the position of the top face of the GIMP in contact and $u_{RB}$ is the bottom position of the rigid body. 

\begin{figure}[ht!]
\centering
\begin{subfigure}[b]{0.49\textwidth}
      \centering
          \includegraphics[width=\textwidth]{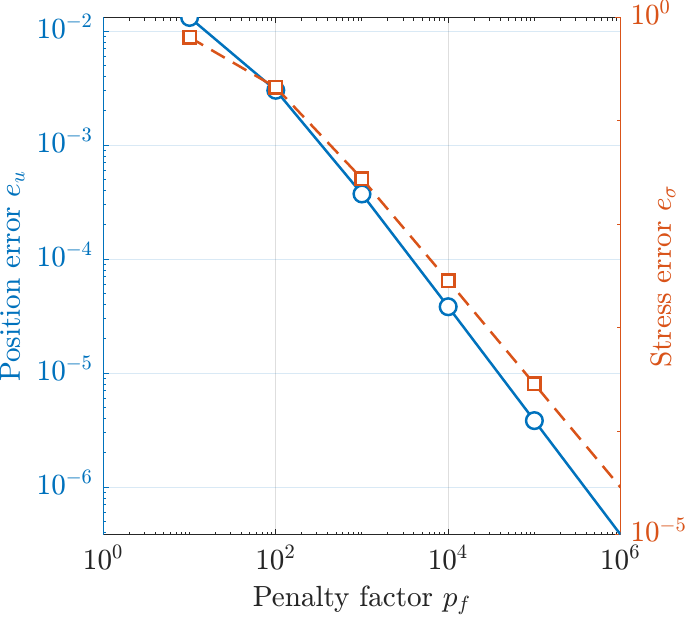}
           \caption{}
         \label{fig:contact column plot error}
          \end{subfigure}
      \hfill
      \begin{subfigure}[b]{0.49\textwidth}
      \centering
          \includegraphics[width=\textwidth]{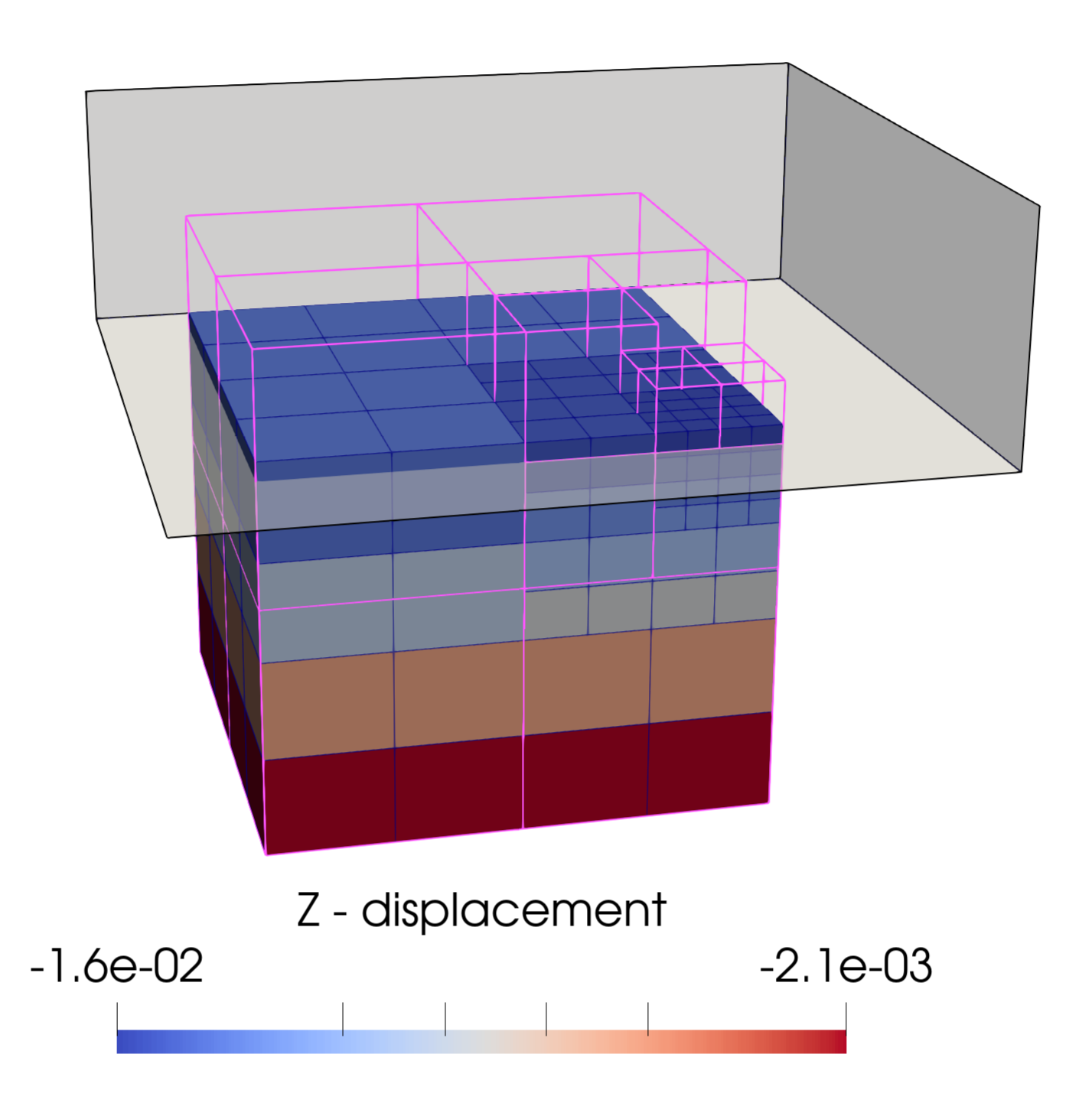}
          \caption{}
          \label{fig:contact column plot stress}
      \end{subfigure}

      \begin{subfigure}[b]{0.4\textwidth}
      \centering
          \includegraphics[width=\textwidth]{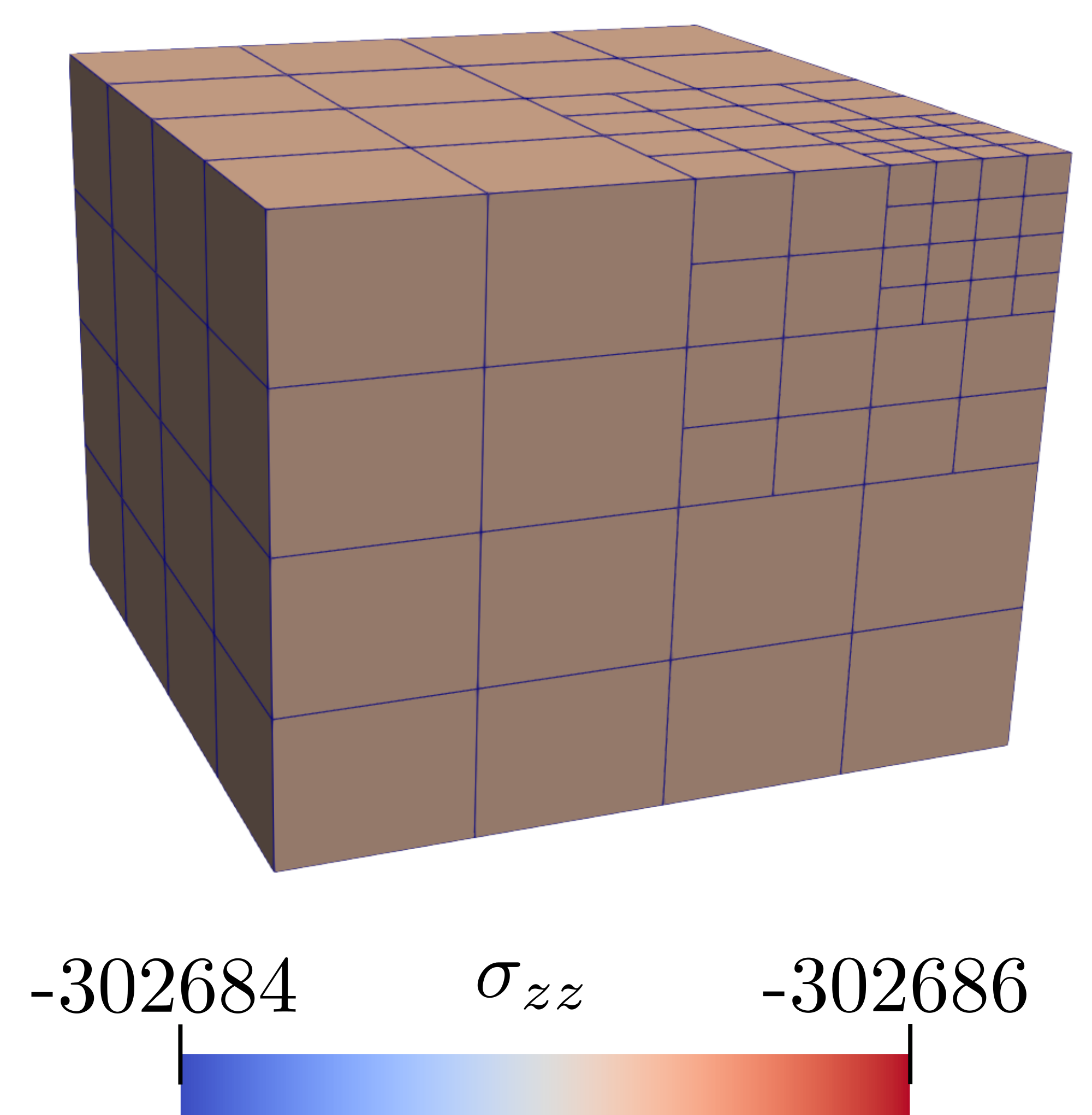}
          \caption{}
          \label{fig:contact column plot stress}
      \end{subfigure}
      \caption{Contact cube: The convergence of the stress and displacement error with $p_f$ is shown in (a) with an example of the final rigid body and GIMP positions, with the mesh shown in pink, provided by (b) and the vertical stress shown in c.}
\end{figure}

The convergence behaviour of these errors is shown in Figure \ref{fig:contact column plot error}, where, as the penalty factor is increased both the error in the overlap across the domain and the stress error decrease with increasing $p_f$.
\ref{fig:column_compression_series_BCs}. 
A graphical example of when $p_f=100$ is shown in Figure \ref{fig:column_compression_series_BCs} with the corresponding vertical stress (Pa) shown in Figure \ref{fig:contact column plot stress}.  The plots clearly show that the position of the top surface and the vertical stress are uniform. This is with GIMPs of varying size, hanging nodes and with the problem experiencing considerable deformation, i.e. a 25$\%$ reduction in height. Figure \ref{fig:column_compression_series_BCs} also demonstrates the overlap that is experienced if $p_f$ is not sufficiently large. This is a particularly stiff problem, being elastic and constrained by its boundary conditions, such that  only when $p_f>1000$ is the error less than $1\%$. However, generally for less constrained problems it is found that $p_f=50$ is generally sufficient.  

\subsection{Sphere rolling down a slope}
\begin{figure}[ht!]
    \centering
    \includegraphics[width=0.75\linewidth]{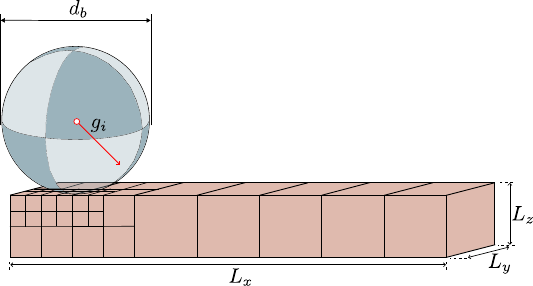}
   \caption{Rolling sphere: Domain setup and dimensions of the rolling sphere problem.}
    \label{fig:Sphere_setup}
\end{figure}

The final problem with an analytical solution is a rigid sphere rolling down an inclined rigid slope. It validates the contact and friction implementation in the presence of an adaptive mesh with hanging nodes. The problem setup, domain, and loading are shown in Figure \ref{fig:Sphere_setup}. The sphere has diameter $d_p = 2.0$ m, mass of $10^4$ kg, and rotational inertia of $4000$ kg~m$^2$.  The sphere is constructed from 3120 triangles arranged on a latitude-longitude grid so that the smallest elements exist at the poles. The sphere's poles are aligned with the plane of rotation and so will contact the GIMPs. The rigid surface over which the sphere rolls is approximated with a stiff domain with a Young's modulus and Poisson's ratio of $10^{9}$ Pa and $0.0$ respectively. A stiff domain is used, rather than fixing all degrees of freedom, to validate the hanging node formulation with contact and friction.  All displacements at the nodes are fixed in $x$ and $y$ with roller boundary conditions applied to all surfaces, except the top surface where there is a homogeneous Neumann boundary condition. The slope has dimensions $L_x = 50$ m and $L_y = L_z = 1$ m, as shown in Figure \ref{fig:Sphere_setup}, the maximum element size is $0.5$ m, whilst the element size at the contact point is set to $0.1$ m. The parameters used for the penalty contact are defined in Section \ref{sec: Rigid body representation and contact approach} but are repeated here for clarity, i.e.
\begin{equation}
    \epsilon_N = 50E_pA_p\quad\text{and}\quad \epsilon_T = 25E_p A_p.
\end{equation}
$\epsilon_N$ and $\epsilon_T$ are the normal and tangential contact penalties, and $E_p$ and $A_p$ are the Young's modulus and contact surface area of the material point $p$ in contact.

To obtain an accurate representation of the slope geometry at a $45^\circ$ incline, the gravitational load is set to $g_i = 9.81 \times [1/\sqrt{2},\, 0,\, -1/\sqrt{2}]^\top$ m/s$^2$. An equivalent setup with an inclined slope and a mesh that adapts as the sphere rolls would misalign the GIMP contact vertices with the correct slope boundary, so here, instead, the gravity vector is tilted.
The analytical solution to this problem is presented as the distance $d_x$ the sphere has moved down the slope as a function of time. The solution has two forms depending if the frictional contact is in a slip or stick state
\begin{equation}
    d_x(t) =
    \begin{cases}
        \dfrac{gt^2}{2}\left[\sin(\theta_s)-\mu\cos(\theta_s)\right] & \text{if slipping, } \tan(\theta_s) > 3.5\mu \\[4pt]
        \dfrac{5gt^2 \sin(\theta_s)}{14} & \text{otherwise (sticking)}
    \end{cases}
\end{equation}
where $g=|g_i|=9.81$ m/s$^2$ is the gravitational acceleration, $\theta_s = 45^\circ$ is the slope angle, $t$ s is time and $\mu\in\{0,0.1,0.2,0.4,1.0\}$ are the friction coefficients considered here.

\begin{figure}
    \centering
    \begin{subfigure}{0.49\linewidth}
        \centering
        \includegraphics[width=\linewidth]{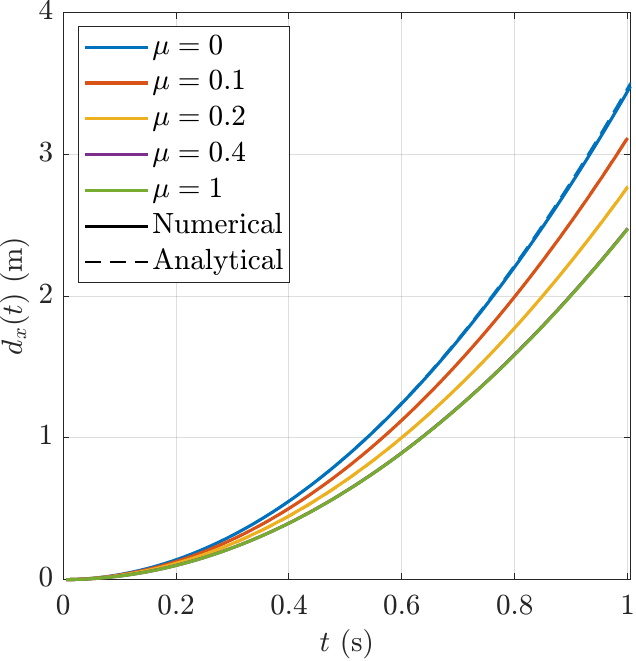}
        \caption{}
        \label{fig:rolling_sphere_results}
    \end{subfigure}
    \hfill
    \begin{subfigure}{0.49\linewidth}
        \centering
        \includegraphics[width=\linewidth]{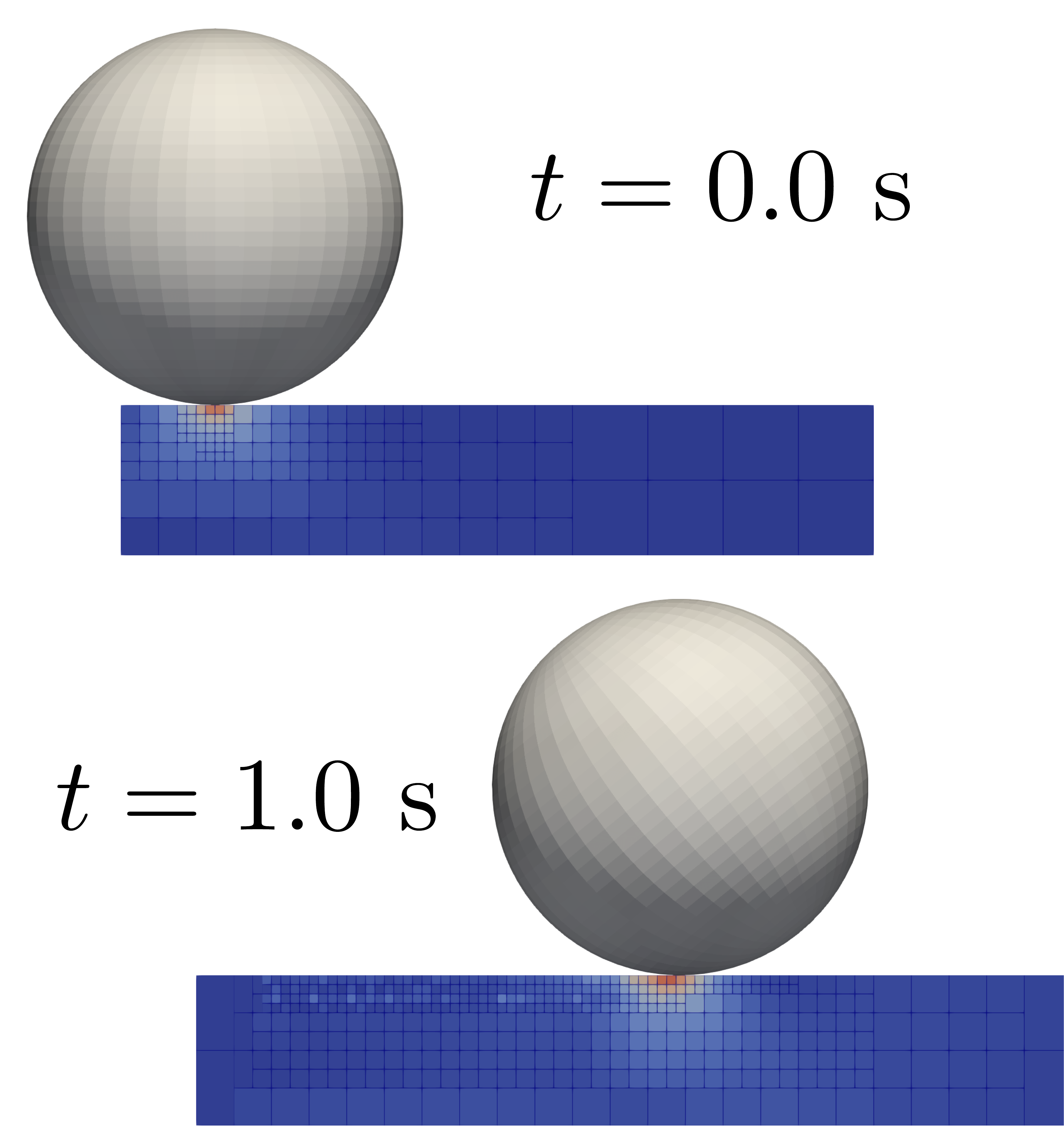}
        \caption{}
        \label{fig:sphere_results_ref}
    \end{subfigure}
    \caption{Rolling sphere: The comparison of simulation results with reference data for different coefficients of friction is shown in (a), (b) shows the sphere and GIMP domains for $\mu=1.0$ where red corresponds to the maximum displacement of $10^{-4}$ m.}
\end{figure}

A comparison of the numerical and analytical results is shown in Figure \ref{fig:rolling_sphere_results} whilst an example of the GIMP distribution, and sphere position is shown in Figure \ref{fig:sphere_results_ref}. The first observation is the excellent agreement between the numerical and analytical results in Figure \ref{fig:rolling_sphere_results}. The results validate the penalty contact and friction GIMP formulation in the presence of hanging nodes for a range of friction coefficients and states. Figure \ref{fig:sphere_results_ref} shows how the GIMPs have been refined as the sphere moves over the domain, initially at $t =0.0$ s, the only GIMPs that are refined are those which are in elements that contain the sphere geometry. As the sphere moves along the slope the refined GIMPs remain. 

\subsection{Cone penetration test}\label{sec: cone penetration test}
Cone Penetration Tests (CPTs) are widely used in geotechnical engineering as an \emph{in-situ} test to determine soil properties. Bird \emph{et al.} \cite{birdanchors2026} showed that a numerical CPT can serve as a first step in predicting other geotechnical problems: material properties are calibrated against experimental CPTs in the soil, and the calibrated parameters are then used on other large deformation soil-structure interaction problems in the knowledge that they provide a good match to the physical CPT response. Modelling CPTs quickly and accurately is therefore a useful tool for material calibration within a broader numerical modelling framework. The aim of this section is to validate that:
\begin{enumerate}
    \item an octree based refinement strategy based on the rigid body location is suitable to achieve accurate results; and
    \item the computational expense is relative small in terms of required hardware and time.
\end{enumerate}

The variation of the cone resistance, $q_c$, with normalised depth $D/r$, where $D$ is the penetration depth and $r$ is the cone radius, is used to validate the numerical results against the experimental data.  Experimental data from \cite{davidson2022physical,Cerfontaine2020} for a dry silica sand obtained from Congleton in the UK with a relative density of $38\%$ is used here. The material is modelled as Hencky hyperelastic-perfectly plastic, with a frictional-cone yield surface using a Willam–Warnke deviatoric section and a non-associated Drucker–Prager flow potential. The numerical properties are populated following the same procedure as in \cite{birdanchors2026,bird_dynamic_2025,BIRD2024106646}. The empirical laws proposed by Brinkgreve \emph{et al.} \cite{brinkgreve2010validation} are used to determine material properties from the relative density. In \cite{birdanchors2026} it was found that the best agreement between numerical and experimental results was obtained when the relative density for the Brinkgreve equations was $32\%$, for an experimental sand of $38\%$. However the numerical CPT results of \cite{bird_dynamic_2025} show that using a value of $38\%$ also produces a reasonably accurate result. A summary of the material properties for the $32\%$ relative density sand is provided in Table \ref{tab: material properties}.
so all looks 
\begin{table}[ht!]
\centering
\caption{Cone penetration test: Material properties.}\label{tab: material properties}
\label{tab:material}
\begin{tabular}{l|c}
\textbf{Property} & \textbf{Value ($R_D = 32\%$)}  \\
\hline
Reference Young's modulus, $E^{ref}_{50}$ (kPa) & 19,200 \\
Density, $\rho$ (kN/m$^3$)                      & 16.3   \\
Poisson's ratio                                 & 0.25   \\
Friction angle ($^\circ$)                       & 32.0   \\
Dilation angle ($^\circ$)                       & 2.0    \\
Apparent cohesion (kPa)                         & 0.3    \\
Coefficient of earth pressure at rest, $K_0$    & 0.47   \\
Stiffness exponent, $m_E$                       & 0.60   \\
\hline
\end{tabular}
\end{table}

The initial state of the material from gravitational loading must also be considered; the Young's modulus therefore varies with confinement pressure using the expression provided by \cite{schanz2019hardening}
\begin{equation}\label{eq: E}
    E_{50} = {E}_{50}^{ref}\left(\frac{{{\sigma}_{v} {K}_{0}}}{{{p}^{ref}}} \right)^{m_E}
    \qquad\text{with}\qquad 
    \sigma_v = d_p^0\rho,
\end{equation}
where $K_0 = 1 - \sin(\phi)$ is the coefficient of earth pressure at rest \cite{jaky1944coefficient}, $\sigma_v$ is the vertical stress, $d_p^0$ is the distance of the GIMP centre from the surface of the sample at the initial state (i.e. the initial depth of the material point, distinct from the cone penetration depth $D$), and $m_E$ is an exponent controlling the variation of stiffness \cite{brinkgreve2010validation}. It is important to note that the $E_{50}$ value at each material point does not change during the simulation, it is fixed by the material point's centre position at the initial state. For the interaction between the cone and the GIMPs the coefficient of friction is set to $0.3$, and normal and tangential contact parameters are set to $\epsilon_N = 50E_pA_p$ and $\epsilon_T = 25E_p A_p$, see Sections \ref{sec: Rigid body representation and contact approach} and \ref{sec: contact problem} for a discussion. The ghost parameters, and the discussion on their calculation, is as  in Section \ref{sec: ghost}.

\begin{figure}[ht!]
    \centering
    \begin{subfigure}[b]{0.35\textwidth}
        \centering
        \includegraphics[width=\textwidth]{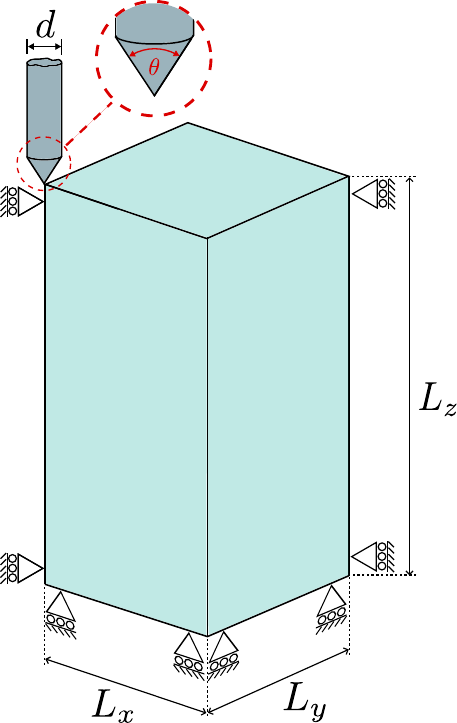}
        \caption{}
       \label{fig:CPT setup}
    \end{subfigure}
    \hspace{0.5cm}
    \begin{subfigure}[b]{0.6\textwidth}
        \centering
        \includegraphics[width=\textwidth]{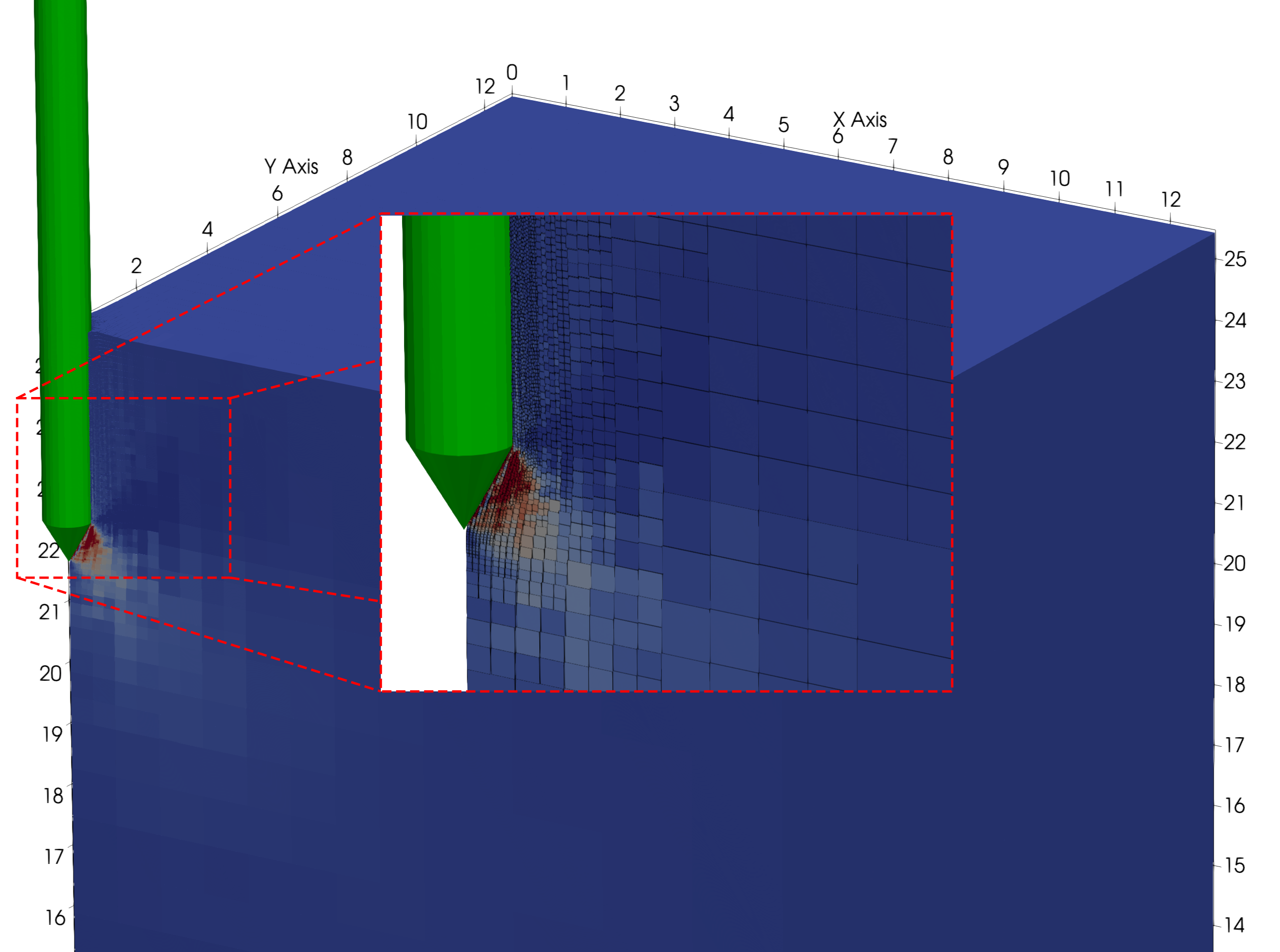}
        \caption{}
        \label{fig:CPT result example}
    \end{subfigure}
    \caption{Cone penetration test: The initial geometry and mesh setup is shown in (a) with the CPT penetrated $1.2$m shown in (b). The displacement magnitude is shown in (b), blue is $0$m and red is $0.5$m.}
\end{figure}

The domain setup for the CPT is provided in Figure \ref{fig:CPT setup}, with dimensions $L_x = L_y = 12.8$~m and $L_z = 24.6$~m, and CPT parameters $r = d/2 = 0.4$~m and $\theta = 60^\circ$, where $d$ is the cone diameter. A coefficient of friction of $\mu = 0.3$ acts between the soil and the CPT, with the numerical penalty contact parameters as defined in Section \ref{sec: Rigid body representation and contact approach}.
The CPT problem has a vertical axis of rotational symmetry through its centre so can be modelled as an axisymmetric problem, as in  Martinelli and Yost \cite{martinelli2022relating,martinelli2022relating,yost2023addressing}, or Bird \emph{et al. }\cite{BIRD2024106646}, or as a 3D slice with tetrahedral elements within the MPM \cite{CECCATO2016440}. With cuboid elements and GIMPs, however, the two-fold symmetry in the $xz$ and $yz$ planes through the CPT centre can be exploited, so only a quarter of the domain need be modelled \cite{bird_dynamic_2025}. Accordingly, roller boundary conditions are imposed on the two symmetry planes and on the remaining lateral and bottom faces. The top face is left as homogeneous Neumann.

In the preceding validation tests the mesh was refined to assess its effect on solution accuracy. Here the adaptive strategy itself is fixed, and the solutions obtained at different minimum element sizes are compared with each other and against the experimental measurements of \cite{davidson2022physical,Cerfontaine2020}. Following the scheme outlined in Section \ref{sec:adaptive stratergies}, elements intersecting the rigid body surface are assigned a size $dx_{min} \in \{0.025, 0.05, 0.1, 0.2\}$~m, and the surrounding region where $dx_{min}$ is set away from the rigid body is $dx_{min}^{region} = 2 \times dx_{min}$.

\begin{figure}[ht!]
    \centering
   \includegraphics[width=1\textwidth]{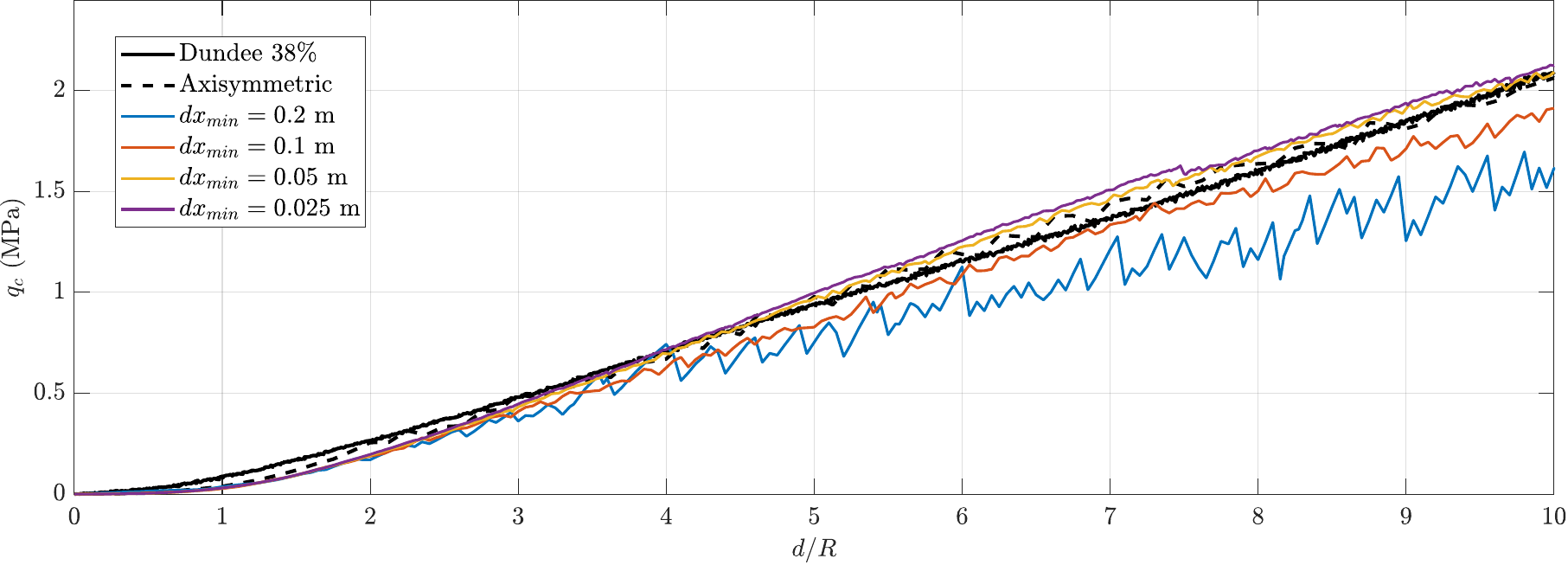}
    \caption{Cone penetration test: Results for CPT tip load with depth.} \label{fig:CPT results}
\end{figure}

The problem is run as a two-stage pseudo static problem. 
\begin{itemize}[leftmargin=2cm]
    \item [Step 1:] In a single load step, apply a gravitational load to the problem. At the end of the step, position the tip of the CPT just above the top of the domain.
    \item [Step 2:] Over 300 load steps vertically displace the CPT a total vertical distance of $-4$ m.
\end{itemize}
The results for the different refinements, and the experimental data of \cite{davidson2022physical,Cerfontaine2020}, are presented in dimensionless form in Figure \ref{fig:CPT results}, with the corresponding simulation times provided in Table \ref{tab:cpt_timing}.

Figure \ref{fig:CPT results} shows two convergence properties: the oscillations in $q_c$ are decreasing and the lines are converging with refinement. However, unlike the other CPT refinement validations, \cite{BIRD2024106646,bird_dynamic_2025},  in general with the GIMPM where it is assumed that a refined mesh will be more compliant, the results here are converging from below suggesting that with refinement the stiffness/strength is increasing. This is caused by the Young's modulus being defined at the GIMP centre point combined with the very large elements that can exist below the CPT. As a result, the centre points of the lowest elements for $dx_{\min}=0.2$ sit higher than those for $dx_{\min}=0.05$. Consequently, for the coarser $dx_{\min}$, the material below the CPT is more compliant; as the mesh is refined, the Young's modulus assigned to these elements becomes larger, so the overall stiffness below the CPT increases. Nonetheless, the results converge with refinement. The wall clock time is presented for each simulation in Table \ref{tab:cpt_timing}, showing that reasonably accurate results, $dx_{\min}=0.1$ m, can be achieved in an hour and highly accurate results, $dx_{\min}=0.05$ m, in approximately $3.5$ hours. For the CPT simulations the linear solve used a direct solver, and it was recorded that the average solve time for $dx_{\min}=0.025$ was approximately $10$ times larger than that for $dx_{\min}=0.1$. Further work would be required to investigate the computation cost with an iterative solver, for example GMRES with a diagonal Jacobi preconditioner.
 
\begin{table}[h]
    \centering
    \caption{Cone penetration test: wall-clock time as a function of background-mesh refinement.}
    \label{tab:cpt_timing}
    \begin{tabular}{c|c}
        {$dx_{\min}$ (m)} & {Wall time (hours:mins)} \\\hline
        0.200 &  0:36   \\
        0.100 &  0:58   \\
        0.050 &  3:26   \\
        0.025 &  32:18 \\
        \hline
    \end{tabular}
\end{table}

\subsection{Anchor penetration}
\label{sec: anchor penetration}

The final validation comprises the modelling of anchor penetration in a level soil bed. Results are compared both to the GIMPM result obtained by Bird \emph{et al.} \cite{birdanchors2026} using a structured mesh, with a $0.1$ m grid size around the anchor, and the experimental data from Sharif \emph{et al.} \cite{sharif}, for an anchor dragged through a sand of $32\%$ relative density. The comparisons expose the variation and sensitivity of the predicted anchor trajectory to the mesh refinement around the contact interface, and the associated computational and environmental (carbon) cost is also investigated. The machine setup is the same as in Section~\ref{sec: cone penetration test}. 

\begin{figure}
    \centering
    \includegraphics[width=0.6\linewidth]{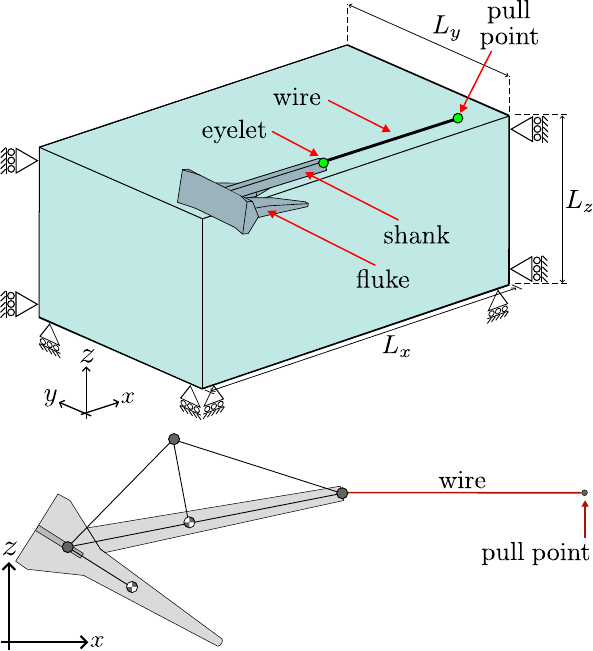}
    \caption{Anchor penetration: the top figure provides the domain dimensions whilst the bottom shows the truss frame used to model the anchor and wire, it also shows where the pull point is applied.}
    \label{fig: anchor_domain}
\end{figure}

The problem domain and anchor setup are provided in Figure~\ref{fig: anchor_domain}. The domain has side lengths $L_y = L_z = 10$ m and $L_x = 100$ m. The material in the domain has the same properties as in Section~\ref{sec: cone penetration test}, with the same variation in Young's modulus with depth, which is also consistent with the material setup in Bird \emph{et al.} \cite{birdanchors2026}. Furthermore, to be consistent with Bird \emph{et al.} \cite{birdanchors2026} the coefficient of friction acting between the anchor and soil is set to 0.3 with normal and tangential contact parameters set to $\epsilon_N = 50E_pA_p$ and $\epsilon_T = 25E_p A_p$ respectively - see Section \ref{sec: Rigid body representation and contact approach}. The stabilising ghost parameters are the same as those defined in Section \ref{sec: ghost}.

This problem is also run in two stages:
\begin{itemize}
    \item[] Step 1: Apply a gravitational load in a single pseudo-static load step;
    \item[] Step 2: Place the anchor on the sand surface, with the pull point at a height of $10$ m. This step is also modelled dynamically and does not progress to step 3 until the vertical oscillations have decreased below the velocity of $10^{-3}$ m/s;  
    \item[] Step 3: The pull point is then moved with a velocity of $(v_x,v_y,v_z)=(0.1,0.0,0.0)$ m/s. This step is modelled dynamically with a time step of $\Delta t=0.01$ s until the anchor has been dragged $19$ m.
\end{itemize}
For all steps the top of the domain has a homogeneous Neumann boundary condition, and for Step 1 all other faces have roller boundary conditions. Although the anchor is dragged $19$ m in Step 2, it is not necessary to model this entire length - instead a partitioned domain method is used, consistent with Bird \emph{et al.} \cite{birdanchors2026}. For the total anchor length $L_a$ m (the difference between the maximum and minimum $x$ coordinates of the rigid body), only $1.5L_a$ m of the domain is modelled in front of the anchor and $0.5L_a$ m behind. The kinematics of the anchor are also modelled using the same method as in \cite{birdanchors2026}, shown at the bottom of Figure~\ref{fig: anchor_domain}, where the anchor comprises two parts hinged together: a shank and a fluke. Their kinematics are modelled with a truss frame in which each truss member has a stiffness of $10^9$ N/m and a nominal nodal mass of $10$ kg.

\begin{table}[ht!]
\centering
\caption{Anchor penetration: total mass, rotational inertia (rot. inertia) and geometric properties of the anchor components, where the CoM position and length are measured from the pivot point.}
\label{tab: anchor mass, intertia and length properties}
\begin{tabular}{r|c|c|c|c}
      & mass (kg)  & rot. inertia (kgm$^2$) & CoM position (m) & length (m) \\ \hline
fluke & 6583.2 & 1100 & 0.131 & 1.7 \\ \hline
shank & 2116.7 & 1350 & 1.272 & 3.3 \\
\end{tabular}
\end{table}
The shank and fluke's total mass, and rotational inertia, properties are provided in Table~\ref{tab: anchor mass, intertia and length properties}. However in this simulation there is a plane of symmetry through the centre of the anchor in the $x$-$z$ plane, and so the mass and rotational inertia are both halved for the analysis, \cite{birdanchors2026}.

\begin{table}[ht!]
    \centering
    \caption{Anchor penetration: list of simulation references, their refinement parameters and their estimated environmental cost; ${^*}$ denotes a value estimated for simulation A, which was stopped early (see text), with the run time extrapolated linearly from the cost per unit drag distance measured over the first $8$ m and the speed up and carbon cost derived from this estimate.}
    \resizebox{\textwidth}{!}{%
    \begin{tabular}{r|c|c|c|c|c|c}
        Simulation ref & $dx_{min}$ & $dx_{min}^{region}$  & run time (hrs) & speed up & carbon cost (kgCO\textsubscript{2}e) & relative carbon cost \\ \hline
        A & 0.1  & 0.3  & 476.9        & 2.5        & 22.38       & 0.104 \\
        B & 0.1  & 0.2  & 218.6        & 5.5        & 10.25       & 0.047 \\
        C & 0.1  & 0.1  & 81.3         & 14.7       & 3.81        & 0.017 \\
        D & 0.2  & 0.6  & 111.0        & 10.8       & 5.20        & 0.024 \\
        E & 0.2  & 0.4  & 71.2         & 16.8       & 3.34        & 0.015 \\
        F & 0.2  & 0.2  & 40.7         & 29.5       & 1.91        & 0.008 \\
    \end{tabular}%
    }
    \label{tab:simulations anchor}
\end{table}

\begin{figure}[ht!]
    \centering
    \includegraphics[width=\linewidth]{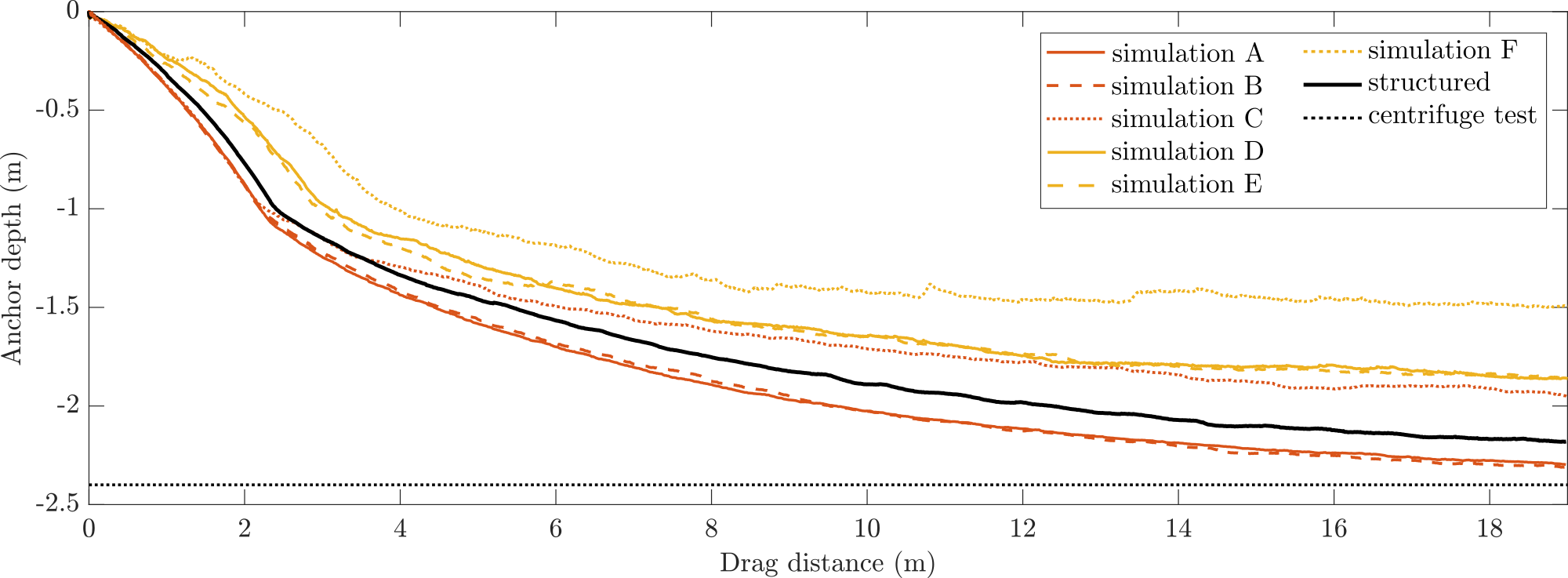}
    \caption{Anchor penetration: comparison of anchor trajectories for different adaptivity schemes.}
    \label{fig: anchor_results}
\end{figure}

\begin{figure}[ht!]
    \centering
    \includegraphics[width=0.9\linewidth]{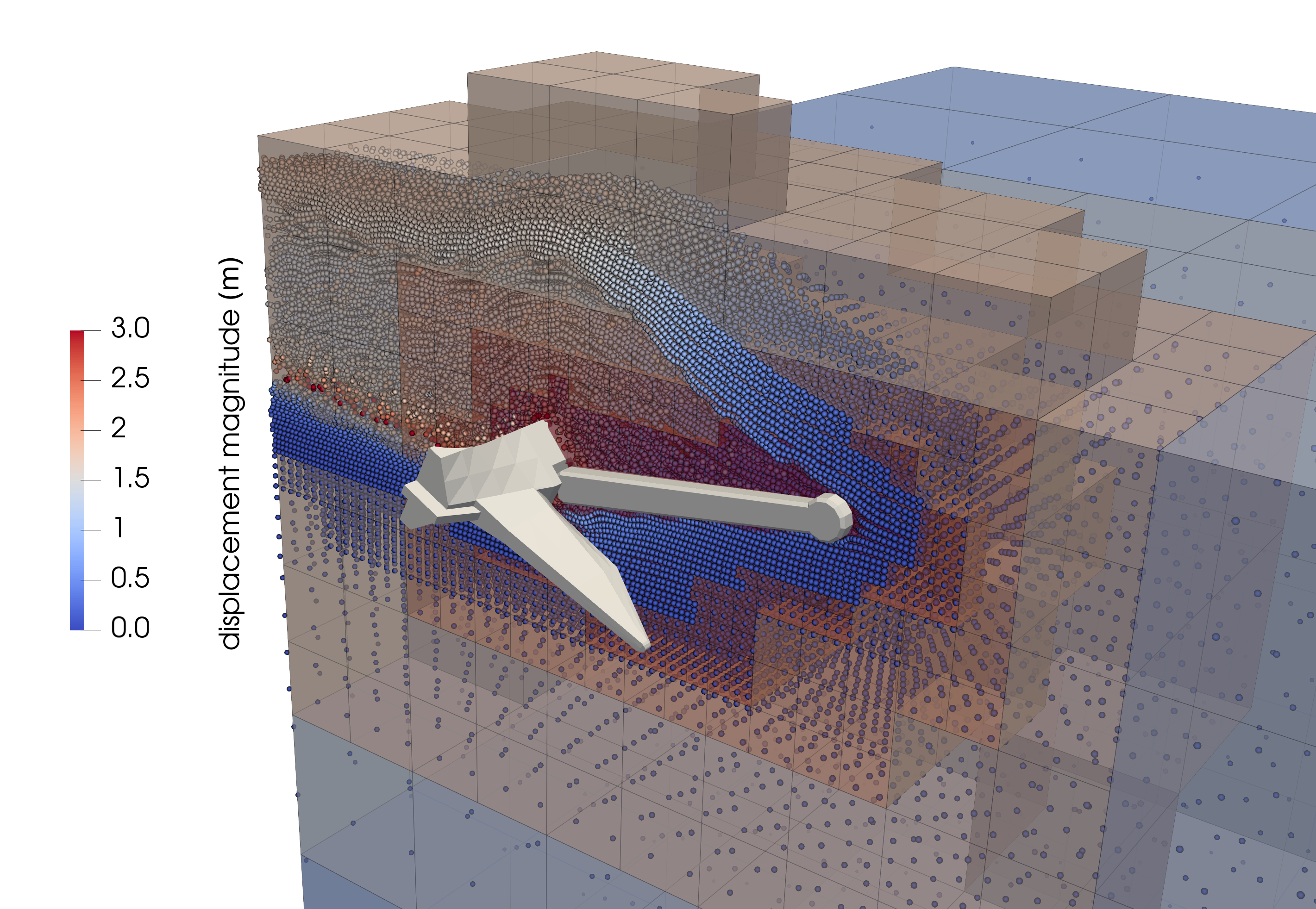}
    \caption{Anchor penetration: Octree background mesh and GIMP distribution for the anchor at a drag distance of $19$ m for Simulation B, where the background mesh is coloured according to refinement age (oldest blue, youngest red).}
    \label{fig: anchor example}
\end{figure}

Table~\ref{tab:simulations anchor} describes the different adaptivity parameters considered here, linked to the six simulations of this anchor penetration problem. Two element sizes were considered, $dx_{min}=\{0.1,0.2\}$, each with two refinement-region sizes $dx_{min}^{region}=\{dx_{min},2dx_{min},3dx_{min}\}$ to demonstrate convergence with a given minimum element size. Table~\ref{tab:simulations anchor} also reports the speed up relative to the structured mesh, which was allocated 40 cores and up to 200 GB of RAM. The structured mesh ran for approximately $50$ days, and the carbon cost reduction is given relative to the structured-mesh cost of $213.95$ kgCO\textsubscript{2}e. The carbon cost was calculated using an online carbon cost calculator \cite{carboncost}. A $100\%$ occupation of the allocated resources for each simulation is assumed when determining the carbon costs, which means that the results in Table~\ref{tab:simulations anchor} are likely to be overestimates. The corresponding anchor trajectories are compared in Figure~\ref{fig: anchor_results} and an example mesh and GIMP distribution around the anchor at a drag distance of $19$ m are shown in Figure~\ref{fig: anchor example}. In Figure~\ref{fig: anchor example}, the centre points of the GIMPs are marked with spheres coloured with their displacement whilst elements in the mesh are shown in the hexahedral form and coloured relative to their age, blue-to-red is old-to-new. The figure shows that as the anchor moves from left-to-right that GIMPs are being refined in front of the anchor and are also associated with small elements, however as the GIMPs move up and behind they stay refined but are now associated with larger elements.

Considering Figure~\ref{fig: anchor_results}, compared with the reference numerical solution, Simulations A and B achieve the same or better agreement with the experimental result than the structured mesh, with Simulation B obtaining this improved accuracy $5.5$ times faster. The octree mesh also gives a substantial reduction in carbon emissions: Simulation B emitted approximately $21$ times less CO\textsubscript{2}e than the structured-mesh analysis, while the approximate solution from Simulation E achieved a $66$-fold reduction. Convergence with respect to the size of the refinement region is also observed: for the simulations with $dx_{min}=0.1$ m, once $dx_{min}^{region} \geq 2,dx_{min}$ there is no further improvement in the predicted trajectory. This result, together with the CPT results (Section~\ref{sec: cone penetration test}), shows that the refinement scheme achieves good accuracy relative to both the structured meshes and the experimental data. The adaptive refinement is not, however, optimal, and an \emph{a posteriori} error estimator would help identify where the mesh is over- and under-refined and this could be the focus of future research. The structured mesh used a local element size of $0.1$ m around the anchor, with a power law of $1.3$ applied to increase the element size away from the anchor in each Cartesian direction, see Bird \emph{et al.} \cite{birdanchors2026} for details. Figure \ref{fig: anchor_results} suggests that the stiffer structured-mesh response, shown by the smaller absolute penetration depth, arises not from the far-field coarsening rate, which is slower for the structured mesh, but from the near-field coarsening, where the structured mesh begins to grow its elements too soon around the anchor. This is supported by the sensitivity of the results to $dx_{min}^{region}$ in Figure \ref{fig: anchor_results}: the predicted trajectory changes significantly between Simulations B$\rightarrow$C and E$\rightarrow$F where the $dx_{min}^{region}$ varies from $2dx_{min}\rightarrow dx_{min}$. 

\section{Conclusion}

This paper has presented the first implicit octree-based adaptive Material Point Method for the large deformation analysis of continuum-rigid body interaction problems. The development was motivated by the high computational and carbon cost of such analyses when modelling practical problems in geotechnics but clearly has wide applicability. The method employs an octree background mesh that automatically adapts both the computational grid and the material point discretisation to the evolving position of the rigid body interacting with the surrounding continuum. Computational resolution, and therefore expense, is focused to the vicinity of the interface; away from the interface the domain is represented with a coarse mesh. This adaptivity reduces the number of degrees of freedom and material points that are used in the analysis, this leads to a reduction in the runtime and overall computational cost. A further advantage of the method is that the user only needs to define the desired minimum element size and the region around the rigid body where this is enforced.

Beyond being the first time that an octree background mesh has been used in an implicit MPM, the two novel methodological contributions that underpin this paper are the discrete calculation of the GIMP basis functions and the extension of ghost stabilisation to non-conforming meshes. As the GIMP basis is computed using discrete integration, the usual piecewise formulation is avoided, which would otherwise become cumbersome and error-prone on a non-conforming mesh. The discrete integration also removes the necessity of building the basis on a hierarchy of meshes that are then coupled together. Instead, only a single mesh is used. As discussed by Coombs \cite{coombs_ghost_2023}, the small-cut issue is a serious drawback of MPM and GIMPM which has been resolved with ghost stabilisation, so here to ensure simulation stability the ghost formulation was extended to non-conforming meshes.

The method was verified and validated on five numerical examples of increasing complexity. The first validation tested the GIMPM's hanging node formulation using a stress convergence study of a column undergoing deformation from self-weight when hanging nodes are present in the mesh. The next validation was for the normal contact formulation, originally presented by \cite{bird_dynamic_2025} for the GIMPM on a conforming mesh, and applied to a problem where hanging nodes were present for the GIMPs  in contact with the rigid body. The next validation extended the contact validation to include friction in the presence of hanging nodes by modelling a sphere rolling down a slope, and was also the first validation to include adaptivity of the mesh driven by the rigid body. The final two validations, a Cone Penetration Test (CPT) and an anchor pull, are compared to experimental data and previous numerical simulations with conforming meshes. The CPT simulation showed excellent agreement with previously run axisymmetric tests and also demonstrated convergence with mesh refinement. A similar observation was made for the anchor simulations, where a more accurate numerical solution was $5.5$ times faster when compared to conforming mesh run. This was also run with reduced machine requirement which overall led to a total carbon reduction of $21$ times.

\section*{Acknowledgements}
This work was supported by the Engineering and Physical Sciences Research Council [grant numbers EP/W000970/1, EP/W000997/1, EP/W000954/1 and UKRI788]. The third author was supported by funding from the Faculty of Science, Durham University. The fourth was supported by the Engineering and Physical Sciences Research Council's Northern Net Zero Accelerator [grant number: EP/Y024052/1]. This work has made use of the Hamilton HPC Service of Durham University. All data created during this research are openly available at \href{https://collections.durham.ac.uk/}{collections.durham.ac.uk} (specific DOI to be confirmed if/when the paper is accepted).

\bibliographystyle{elsarticle-num} 
\bibliography{cas-refs.bib}
\end{document}